\documentclass[11pt]{amsart} \usepackage[all]{xy} \usepackage{enumitem} \usepackage{mathtools} \usepackage{mathrsfs} \usepackage{amssymb} \usepackage{amsmath,amsfonts,amsthm} \usepackage{graphicx} \usepackage{upgreek} \usepackage[toc,page]{appendix} \usepackage{bbm} \usepackage{esint} \usepackage{scalerel} \usepackage{stackengine,wasysym} \usepackage{verbatim}  \usepackage{hyperref} \usepackage[dvipsnames]{xcolor} \hypersetup{colorlinks=true, linktocpage=true,citecolor=violet, linkcolor=violet, urlcolor=violet} \usepackage{geometry} 
\usepackage{tikz}
\usepackage{stmaryrd}
\newtheorem{teo}{Theorem}[section]
\newtheorem{defi}[teo]{Definition}
\newtheorem{prop}[teo]{Proposition}
\newtheorem{cor}[teo]{Corollary}
\newtheorem{lemma}[teo]{Lemma}

\newtheorem{conjecture}[teo]{Conjecture}
\newtheorem{rmk}[teo]{Remark}
\newtheorem{mainteo}{Theorem}

\newcommand{\R}{\mathbb{R}}
\newcommand{\Ric}{\text{Ric}}
\newcommand{\vol}{\text{Vol}}
\usepackage[english]{babel}
\usepackage{ragged2e}
\usepackage{amsmath}
\usepackage{mathtools}
\usepackage{xcolor}
\usepackage{amsfonts}
\usepackage{amsthm}
\usepackage{thmtools}
\usepackage{float}
\usetikzlibrary{arrows.meta}
\definecolor{mylavender}{RGB}{204,204,251}
\def\TA{1.6}   
\def\TB{0.85}  
\def\SC{3.05}  
\def\SR{0.8}   

\colorlet{myfill}{blue!20}

\newcommand{\TorusHolePath}{%
  (-0.62,0.06)
  .. controls (-0.3,-0.22) and (0.3,-0.22) ..
  (0.62,0.06)
  -- (0.47,-0.04)
  .. controls (0.2,0.13) and (-0.2,0.13) ..
  (-0.47,-0.04)
  -- cycle
}

\newcommand{\TorusHole}{%
  \draw (-0.62,0.06)
    .. controls (-0.3,-0.22) and (0.3,-0.22) ..
    (0.62,0.06);

  \draw (-0.47,-0.04)
    .. controls (-0.2,0.13) and (0.2,0.13) ..
    (0.47,-0.04);
}

\newcommand{\TorusNeckSphere}[2]{%
  \begin{scope}[xshift=#1cm]
    \pgfmathsetmacro{\h}{#2}
    \pgfmathsetmacro{\t}{asin(\h/\TB)}
    \pgfmathsetmacro{\s}{asin(\h/\SR)}
    \pgfmathsetmacro{\xa}{\TA*cos(\t)}
    \pgfmathsetmacro{\xb}{\SC-\SR*cos(\s)}
    \pgfmathsetmacro{\d}{0.3*(\xb-\xa)}
    \pgfmathsetmacro{\hm}{0.62*\h}

    \path[fill=myfill, even odd rule]
      (\xa,\h)
      arc[start angle=\t, end angle=360-\t,
          x radius=\TA, y radius=\TB]
      .. controls (\xa+\d,-\hm)
          and (\xb-\d,-\hm) .. (\xb,-\h)
      arc[start angle=180+\s, end angle=540-\s,
          radius=\SR]
      .. controls (\xb-\d,\hm)
          and (\xa+\d,\hm) .. (\xa,\h)
      -- cycle
      \TorusHolePath;

    \draw
      (\xa,\h)
      arc[start angle=\t, end angle=360-\t,
          x radius=\TA, y radius=\TB]
      .. controls (\xa+\d,-\hm)
          and (\xb-\d,-\hm) .. (\xb,-\h)
      arc[start angle=180+\s, end angle=540-\s,
          radius=\SR]
      .. controls (\xb-\d,\hm)
          and (\xa+\d,\hm) .. (\xa,\h)
      -- cycle;

    \draw[densely dashed]
      (\xa,\h)
      arc[start angle=90, end angle=270,
          x radius=0.18*\h, y radius=\h];

    \draw[densely dashed]
      (\xb,-\h)
      arc[start angle=-90, end angle=90,
          x radius=0.18*\h, y radius=\h];

    \TorusHole
  \end{scope}
}

\newcommand{\TorusSphereApart}[1]{%
  \begin{scope}[xshift=#1cm]

    \path[fill=myfill, even odd rule]
      (0,0) ellipse[x radius=\TA, y radius=\TB]
      \TorusHolePath;

    \draw (0,0)
      ellipse[x radius=\TA, y radius=\TB];

    \filldraw[fill=myfill, draw=black]
      (\SC+0.25,0) circle[radius=\SR];

    \TorusHole

  \end{scope}
}

\def\SphR{0.8}   
\def\SphC{2.6}   

\colorlet{myfill}{blue!20}

\newcommand{\SphereNeckSphere}[2]{%
  \begin{scope}[xshift=#1cm]
    \pgfmathsetmacro{\h}{#2}
    \pgfmathsetmacro{\s}{asin(\h/\SphR)}
    \pgfmathsetmacro{\xa}{\SphR*cos(\s)}
    \pgfmathsetmacro{\xb}{\SphC-\SphR*cos(\s)}
    \pgfmathsetmacro{\d}{0.3*(\xb-\xa)}
    \pgfmathsetmacro{\hm}{0.62*\h}

    \filldraw[fill=myfill, draw=black]
      (\xa,\h)
      arc[start angle=\s, end angle=360-\s, radius=\SphR]
      .. controls (\xa+\d,-\hm) and (\xb-\d,-\hm) .. (\xb,-\h)
      arc[start angle=180+\s, end angle=540-\s, radius=\SphR]
      .. controls (\xb-\d,\hm) and (\xa+\d,\hm) .. (\xa,\h)
      -- cycle;

    \draw[densely dashed] (\xa,\h)  arc[start angle=90,  end angle=270,
                                        x radius=0.18*\h, y radius=\h];
    \draw[densely dashed] (\xb,-\h) arc[start angle=-90, end angle=90,
                                        x radius=0.18*\h, y radius=\h];
  \end{scope}
}

\newcommand{\SpheresSegment}[2]{%
  \begin{scope}[xshift=#1cm, yshift=#2cm]
    \filldraw[fill=myfill, draw=black] (0,0) circle[radius=\SphR];
    \filldraw[fill=myfill, draw=black] (\SphC,0) circle[radius=\SphR];
    \draw (\SphR,0) -- (\SphC-\SphR,0);
  \end{scope}
}

\newcommand{\SpheresApart}[2]{%
  \begin{scope}[xshift=#1cm, yshift=#2cm]
    \filldraw[fill=myfill, draw=black] (0,0) circle[radius=\SphR];
    \filldraw[fill=myfill, draw=black] (\SphC,0) circle[radius=\SphR];
  \end{scope}
}

\title{Stability for the Bray-Brendle-Neves Scalar Sphere Rigidity Theorem}
\author{Carolina Ruviaro}
\address{IMPA - Instituto de Matemática Pura e Aplicada, Rio de Janeiro, RJ, Brasil.} \email{carolinaruviaro@gmail.com} 

\usepackage{setspace}
\begin{document}
\raggedbottom

\begin{abstract} 
Let $(M^3,g)$ be a closed and oriented 3-dimensional Riemannian manifold with $\pi_2(M)\neq 0$ and let $\mathcal{A}_{\pi_2}(M,g)$ denote the least area among homotopically nontrivial two-spheres in $M$. By a theorem of Bray-Brendle-Neves, $\mathcal{A}_{\pi_2}(M,g)\inf_M R_g \leq 8\pi$ with equality characterizing the cylindrical geometry $S^2\times\mathbb{R}$. Motivated by Sormani's questions on volume preserving intrinsic flat stability for scalar curvature rigidity problems, we study the almost equality regime of this theorem. We build, using Sweeney's tunnel construction, a sequence satisfying $\mathcal{A}_{\pi_2}(M_j,g_j)\inf_{M_j} R_{g_j} \rightarrow 8\pi$ with uniformly bounded volume and diameter whose volume preserving intrinsic flat limit is not the Bray-Brendle-Neves rigidity space, showing that stability fails in general. We then prove stability in two special cases: for sequences of rotationally symmetric metrics (with non-collapsing assumptions) and for sequences of manifolds with nonnegative Ricci curvature.
\end{abstract} 

\maketitle \setcounter{tocdepth}{1} 

\section{Introduction}

Rigidity theorems are among the most interesting results in Riemannian Geometry. A typical rigidity theorem asserts that equality in a sharp geometric  inequality forces the underlying manifold to be isometric, or closely related, to a distinguished model space. Once such a theorem is understood, a natural subsequent problem is to understand  whether such a conclusion persists under almost equality. More precisely, one asks whether a sequence of manifolds for which the relevant geometric quantity approaches its extremal value must converge, in a suitable sense, to one of the corresponding rigidity models.

In the context of scalar curvature, stability questions are particularly delicate. The scalar curvature of a Riemannian manifold is a pointwise trace quantity and places much weaker restrictions on the global geometry than, for example, sectional or Ricci curvature. In particular, positive scalar curvature is compatible with thin tunnels, wells, bubbles  and other local modifications which can have dramatic effects on metric convergence while occupying arbitrarily small volume. For this reason, one cannot in general expect Gromov-Hausdorff compactness under scalar curvature bounds alone, and other notions of convergence, such as the intrinsic flat introduced by Sormani and Wenger \cite{SormaniWenger2011}, have become important tools in the study of these problems.

The purpose of this paper is to study stability for the scalar sphere rigidity theorem of Bray-Brendle-Neves. Let $(M^3,g)$ be a closed orientable Riemannian $3$-manifold with $\pi_2(M)\neq 0$. Denote $\mathcal{F}$ the set of all smooth maps $f \colon S^2 \rightarrow M$ that correspond to a nonzero class of $\pi_2(M)$ and define \[\mathcal{A}_{\pi_2}(M,g) \coloneq \inf \{\text{Area}_{f^*g}(S^2) \colon f \in \mathcal{F}\}.\] In this context they proved the following.

\begin{teo}[Theorem 1 in \cite{BrayBrendleNeves2010}] \label{BBN theorem}
    If $R_g$ denotes the scalar curvature of $(M^3,g)$, we have \[\mathcal A_{\pi_2}(M,g)\inf_M R_g\leq 8\pi.\] Moreover, equality holds if and only if the universal cover of $(M^3,g)$ is isometric to the standard cylinder $S^2 \times \R$ up to scaling.
\end{teo}

Thus the quantity $\mathcal{A}_{\pi_2}(M,g)\inf_M R_g$ detects a global splitting phenomenon at the extremal value $8\pi$. It is therefore natural to ask whether the Bray-Brendle-Neves theorem admits an almost rigidity version. Namely, suppose that a sequence of closed and oriented $3$-dimensional Riemannian manifolds $(M_j,g_j)$ satisfies \[\mathcal{A}_{\pi_2}(M_j,g_j)\inf_{M_j}R_{g_j}\xrightarrow{\ j\rightarrow +\infty \ } 8\pi.\] We can ask ourselves: 
\begin{center}
\textit{under which additional hypotheses must a subsequence converge to a space associated with the equality case of Theorem \ref{BBN theorem}? Which notion of convergence should be expected?}
\end{center}
This question was also raised by Sormani in the setting of intrinsic flat convergence. In \cite{Sormani2023Conjectures}, she asked whether Theorem \ref{BBN theorem} admits a volume preserving intrinsic flat stability formulation. More generally, this problem is part of a broader effort to understand what additional geometric assumptions are needed to avoid degeneration phenomena that can occur under lower bounds on scalar curvature. Natural conditions in this setting include lower bounds for the area of minimal surfaces, non-collapsing assumptions and diameter bounds.

The first result of this paper shows that almost equality in Theorem \ref{BBN theorem}, even together with uniform global volume and diameter bounds, is not sufficient to guarantee stability. The counterexample construction is based on the tunnel technique developed by Sweeney \cite{Sweeney2025}, which is an improvement of Gromov-Lawson's construction \cite{GromovLawson1980} (see also Schoen-Yau \cite{SchoenYau1979}), and on a technical theorem that extends Lakzian-Sormani \cite[Theorem 4.6]{LakzianSormani2013} in the context of disconnected Riemannian manifolds, which may have independent interest.

\begin{restatable}{mainteo}{theoremA}\label{tunnel counterexample}
There exists a sequence of Riemannian manifolds $(M_j,g_j)$ and positive constants $v_0,V,D$ independent of $j \in \mathbb N$ such that \[\mathcal{A}_{\pi_2}(M_j,g_j)\inf_{M_j} R_{g_j} \xrightarrow{\ j \rightarrow +\infty \ } 8\pi, \quad v_0 \leq \vol(M_j) \leq V, \quad \text{diam}(M_j) \leq D\] and $M_j$ converges in the volume preserving intrinsic flat sense to the disjoint union of \ $S^2\times S^1$ and $S^3$. 
\end{restatable}

We begin with a sequence of manifolds $X_j=S^2\times S^1_\ell$ endowed with the product metric $h_j \coloneqq r_j^2 g_{S^2}+dt^2$, where $g_{S^2}$ denotes the round metric, $\ell>0$ is the $S^1$-length and the radii $r_j\to 1$. One has \[R_{h_j}=\frac{2}{r_j^2} \qquad \text{and} \qquad \mathcal{A}_{\pi_2}(X_j,h_j)=4\pi r_j^2.\] We then attach to $X_j$ a round $3$-sphere through a tunnel whose radius tends to zero. Sweeney's construction (see Proposition \ref{tunnel} below) allows the tunnel to be chosen with arbitrarily small volume and with an arbitrarily small loss in the scalar curvature lower bound. At the same time, a calibration argument shows that the least area of a nontrivial spherical homotopy class remains close to that of the original $S^2\times S^1$ component. Therefore the resulting manifolds $(M_j,g_j)$ satisfy the almost rigidity condition \[ \mathcal{A}_{\pi_2}(M_j,g_j) \inf_{M_j}R_{g_j} \xrightarrow{\ j \rightarrow +\infty \ } 8\pi.\] 

Nevertheless the additional component $S^3$ has positive volume independent of $j$ and hence it survives in the volume preserving intrinsic flat limit. The limiting space therefore contains both the $S^2 \times S^1$ component and the bubble. To prove this convergence, we make use of Theorem \ref{disjoint Lakzian-Sormani} below, where, in a broader setting, it is established a way to estimate the intrinsic flat distance between two manifolds joined by a tunnel. This theorem extends Lakzian-Sormani \cite[Theorem 4.6]{LakzianSormani2013} using methods of Sormani-Wenger \cite[Section 5]{SormaniWenger2011} and may be of particular interest. To apply Theorem \ref{disjoint Lakzian-Sormani} to our construction, we need to consider two Riemannian metrics on each manifold. We therefore state the theorem in this more general situation. However, the reader may take the two metrics to be equal, which gives a simpler version of the result.

\begin{mainteo}\label{disjoint Lakzian-Sormani}
Consider $N_1$ and $N_2$ two connected, closed and oriented $n$-manifolds, $n\geq 2$. Let $h_i$ and $\hat{h}_i$ be two smooth Riemannian metrics on $N_i$ for $i\in \{1,2\}$. Fix two points $p_i\in N_i$ and let $B^{\hat{h}_i}(p_i)\subset N_i$ be a small open $\hat{h}_i$-geodesic ball centered at $p_i$ with smooth boundary. Set \[ U_i \coloneqq N_i\setminus B^{\hat{h}_i}(p_i) \qquad \text{and} \qquad U\coloneqq U_1\sqcup U_2.\] Suppose that $(M,g)$ is a closed and oriented Riemannian $n$-manifold obtained by joining $U_1$ and $U_2$ through a smooth Riemannian tunnel diffeomorphic to $S^{n-1}\times[0,1]$. Assume that the natural identification $F\colon U\longrightarrow M$ is a smooth orientation-preserving embedding satisfying $(F\vert_{U_i})^*g = \hat{h}_i\vert_{U_i}$ for $i \in \{1,2\}$.

Fix $d>0$ and let $(Y,d_Y)$ be the length space obtained by joining $(N_1,h_1)$ and $(N_2,h_2)$ by a segment of length $d$, attached at $p_1$ and $p_2$, and equipped with the induced length metric. Define \[ X\coloneqq N_1\sqcup N_2, \qquad d_X\coloneqq d_Y\vert_{X\times X} \] and consider the $n$-dimensional integral currents \[ T_M\coloneqq\llbracket M\rrbracket \qquad \text{and} \qquad T_X\coloneqq \llbracket N_1\rrbracket + \llbracket N_2\rrbracket \] with compatible orientations. For each $i\in\{1,2\}$, let $L_i\geq 1$ be a constant such that $h_i \leq L_i^2 \hat{h}_i$. Setting \[ \lambda\coloneqq
\sup_{x,y\in U} \left\vert d_g(F(x),F(y)) - d_X(x,y) \right\vert,\] the intrinsic flat distance between $(M, d_g, T_M)$ and $(X, d_X, T_X)$ satisfies the following estimate
\begin{align*}
d_{\mathcal F} \bigl((M,d_g,&\ T_M), (X, d_X,T_X)\bigr) \leq \frac{n+1}{2}\lambda \left(L_1^n \vol_{\hat{h}_1}(U_1) + L_2^n\vol_{\hat{h}_2}(U_2) \right)\\
& + \dfrac{n}{2}\lambda \left( L_1^{n-1} \text{Area}_{\hat{h}_1}(\partial U_1) + L_2^{n-1} \text{Area}_{\hat{h}_2}(\partial U_2) \right) + \vol_g(M \setminus F(U)) + \mathcal H^n_{d_X}(X \setminus U).
\end{align*}
\end{mainteo}

An important feature of the preceding construction in Theorem \ref{tunnel counterexample} is that the topology may be kept fixed. Indeed, by attaching a $3$-sphere we observe that $M_j\cong (S^2\times S^1)\# S^3\cong S^2\times S^1.$ Thus the failure of stability is not caused by a change in topology. This phenomenon, however, has already been identified by Xu \cite{Xu2023TopologicalGap} who provided a topological stability version of Theorem \ref{BBN theorem} (see Theorem \ref{Theorem of Xu} in Section \ref{Section Ric}). 

Therefore Theorem \ref{tunnel counterexample} gives an answer to the following question posed by Sormani:

\begin{conjecture}[Question 11.6 in \cite{Sormani2023Conjectures}]\label{Sormani's Conjecture}
    Can we formulate and prove a volume intrinsic flat geometric stability conjecture for the Cover Splitting Scalar Rigidity Theorem of Bray, Brendle and Neves \cite{BrayBrendleNeves2010}?
\end{conjecture}

Once our construction shows that it is not possible to have a positive answer to this question in its most direct form, we wonder what additional hypotheses could make it true. We first show that for rotationally symmetric metrics the answer to the question is yes under the volume preserving intrinsic flat and also Gromov-Hausdorff convergence. That is, we consider metrics on $S^2\times S^1$ of the form \[g_{\rho} \coloneqq \rho(t)^2g_{S^2}+dt^2\] where $\rho$ is a positive function on $S^1$. When it is necessary to specify the length of the $S^1$-factor we denote by $S^2 \times S^1_\ell$, $\ell>0$. 

\begin{mainteo}\label{Theorem C}
    Let $M_j = S^2 \times S^1_{\ell_j}$ be a sequence of closed and oriented manifolds with rotationally symmetric metric $g_{\rho_j}$ such that:
    \begin{enumerate}
        \item[1)] diam($M_j$) $\leq \ L$  \ \ and \ \  $\vol(M_j) \geq v>0$;
        \item[2)] $\inf_{M_j} R_{g_{\rho_j}} =2$;
        \item[3)] $\mathcal{A}_{\pi_2}(M_j, g_{\rho_j}) \geq 4\pi(1-\varepsilon_j)$ with $\varepsilon_j \rightarrow 0$,
    \end{enumerate}
    where $L$ and $v$ are positive constants that do not depend on $j$, $\varepsilon_j \in (0,1)$ and $\ell_j>0$. Then there exists a subsequence $(M_{j_k})_k$ that converges to the standard $S^2\times S^1_{\ell_\infty}$, $\ell_\infty \neq 0$, in the sense of Gromov-Hausdorff and also in the volume preserving intrinsic flat sense.
\end{mainteo}

The mechanism behind this result is different from simply fixing the topology. What prevents bubbling in the rotationally symmetric setting is instead the rigidity imposed by the warped product ansatz: the entire geometry is encoded by the single function $\rho_j$ and the almost equality condition controls this function globally.

Before presenting the second positive result, we quickly discuss the significance of the non-collapsing assumptions that have appeared. Even within the model family $S^2\times S^1_{\ell_j}$ with canonical metric, two elementary degenerations are naturally found exploring the role of $\ell_j$. For $S^1$-lengths approaching zero, we see that \[\vol(S^2\times S^1_{\ell_j})\xrightarrow{\ \ell_j \longrightarrow 0\ } 0\] and the sequence collapses in the Gromov-Hausdorff sense to $S^2$. On the other hand, if $\ell_j\to\infty$, the diameters diverge and no compact Gromov-Hausdorff limit exists. Thus lower volume bounds and upper diameter bounds are natural assumptions if one wants a compact $3$-dimensional limit. Theorem \ref{tunnel counterexample} shows, however, that these assumptions in the general case do not prevent bubbling and therefore they are not sufficient by themselves.

Another positive answer proved in this paper is for the case where it is assumed $\Ric \geq 0$. In this situation the notion of convergence is the pointed Gromov-Hausdorff of the universal cover.

\begin{restatable}{mainteo}{theoremC}\label{Ricci theorem}
Let $(M_j,g_j)_j$ be a sequence of closed and oriented Riemannian $3$-manifolds with $\pi_2(M_j) \neq 0$ for all $j \in \mathbb N$. Suppose that: 
    \begin{enumerate}
        \item[1)] $\inf_{M_j} R_{g_j} = 2$;
        \item [2)] $\mathcal{A}_{\pi_2}(M_j,g_j) \geq 4\pi(1-\varepsilon_j)$ with $\varepsilon_j \rightarrow0$;
        \item[3)] $\text{Ric}_{g_j}\geq 0$.
    \end{enumerate}
    Then, up to a subsequence, for all $p \in S^2\times \R$ there exists a sequence of points $p_j \in \widetilde{M}_j$ in the universal cover of $M_j$ such that \[(\widetilde{M}_j,p_j) \xrightarrow{ \ pGH \ } (S^2 \times \R, p)\] where $S^2 \times \R$ is the cylinder with canonical product metric $g_{S^2}+dt^2$.
\end{restatable}

The nonnegative assumption on Ricci curvature along with the already known topology for the universal covers $\widetilde{M}_j$ of the manifolds given by Xu \cite{Xu2023TopologicalGap} allow us to invoke Cheeger-Gromoll splitting theorem. We show, using the theory of Alexandrov spaces, that the $S^2$-factor of these covers converges to the round $2$-sphere. Since $\widetilde{M}_j$ has a product metric, this allows us to conclude, under the pointed Gromov-Hausdorff convergence, that the covers have as limit the canonical cylinder.

The paper is organized as follows. In Section \ref{Section Counterexample} we prove  Theorem \ref{disjoint Lakzian-Sormani} and, with this, we are able to construct the counterexample to Conjecture \ref{Sormani's Conjecture}, thus proving Theorem \ref{tunnel counterexample}. In Section \ref{Section warped} we prove that in the rotationally symmetric case we have Gromov-Hausdorff and volume preserving intrinsic flat convergence, verifying Theorem \ref{Theorem C}. This is done in three parts: first, since the geometry of these manifolds is encoded by a warping function, we give a $C^1$ estimate for these functions. This allows to prove, in Subsection \ref{Subsection GH}, the Gromov-Hausdorff convergence of a subsequence. In Subsection \ref{Subsection IF} we note that the geometric assumptions guarantee a subsequence converging in the intrinsic flat sense. But in this case the Ricci curvature is uniformly bounded below and so the limit in this notion of convergence will coincide with the Gromov-Hausdorff one, solving the problem for rotationally symmetric metrics. Finally, Section \ref{Section Ric} treats the  case of nonnegative Ricci curvature, where it is proved Theorem \ref{Ricci theorem}, which gives an affirmative answer for the stability question in the pointed Gromov-Hausdorff case. For readers not familiar with the concepts of convergence treated in this article, in Appendix \ref{Appendix GH} we recall the definition and basic properties of Gromov-Hausdorff convergence, and in Appendix \ref{Appendix IF} we review intrinsic flat and volume preserving intrinsic flat convergence.

\section*{Acknowledgements}

The author would like to thank her advisor, Lucas Ambrozio, for numerous helpful discussions concerning this problem. The author is also grateful to Christina Sormani for her helpful suggestions and explanations regarding several concepts discussed in this work. This study was financed in part by the Coordenação de Aperfeiçoamento de Pessoal de Nível Superior - Brasil (CAPES) - Finance Code 001 and FAPERJ - Fundação Carlos Chagas Filho de Amparo à Pesquisa do Estado do Rio de Janeiro, Processo SEI-260003/001431/2026.

\section{The Counterexample}\label{Section Counterexample}
The first ingredient in the construction of the counterexample to Conjecture \ref{Sormani's Conjecture} is the tunnel technique developed by Sweeney \cite{Sweeney2025}, which was inspired by Gromov-Lawson's construction \cite{GromovLawson1980}. This construction allows us to connect two Riemannian manifolds by removing sufficiently small geodesic balls and inserting a thin neck while controlling the geometry of the resulting metric. The tunnel can be chosen so that the scalar curvature lower bound decreases by an arbitrarily small amount whereas both the volume of the cylindrical region inserted and the size of the regions modified around the gluing points tend to zero. 

\begin{prop}[Proposition 3.2 in \cite{Sweeney2025}] \label{tunnel}
Let $(M^n,g)$, $n\geq 3$, be a (possibly disconnected) complete Riemannian manifold with scalar curvature $R_g \geq \kappa$. Let $\delta>0$ be small enough, $j\in\mathbb{N}$ and $d\geq 0$. Fix two balls $B_g(p,2\delta)$ and $B_g(p',2\delta)$ in $(M^n,g)$. Then we can construct a new complete Riemannian manifold $P^n$, where we remove the two balls and glue in a cylindrical region $(T_j,g_j)$ diffeomorphic to $S^{n-1}\times[0,1]$, \[P^n = M^n\setminus \left(B_g(p,2\delta)\cup B_g(p',2\delta)\right) \cup T_j.\] Furthermore, the following properties are satisfied:
\begin{enumerate}
    \item $P$ has scalar curvature $R_P>\kappa-\frac{1}{j};$
    \item $g_j\vert _E=g\vert _E$ and $ g_j\vert _{E'}=g\vert _{E'}$ where $E=B_g(p,2\delta)\setminus B_g(p,\delta)$ and $E' = B_g(p',2\delta)\setminus B_g(p',\delta)$ are identified with subsets of $P$;
    \item There exists a constant $C>0$, independent of $j$ and $d$, such that \[d \leq \text{diam}(T_j) < C(\delta+d) \quad\text{and} \quad \vol(T_j) < C\left(\delta^n+d\delta^{n-1}\right);\]
    \item  Outside the tunnel the metric agrees with the original metric $g$ under the natural identification of the corresponding regions.
\end{enumerate}
\end{prop}

\begin{figure}[H]
  \centering
  \resizebox{1\linewidth}{!}{%
  \begin{tikzpicture}[
    line width=0.6pt,
    >={Stealth[length=2.5mm,width=1.5mm]}
]
    \TorusNeckSphere{0}{0.40}
    \draw[->] (4.1,0) -- (4.7,0);
    \TorusNeckSphere{6.4}{0.28}
    \draw[->] (10.5,0) -- (11.1,0);
    \TorusNeckSphere{12.8}{0.12}
    \draw[->] (16.9,0) -- (17.4,0);
    \node at (17.85,0) {\LARGE$\cdots$};
    \draw[->] (18.3,0) -- (18.8,0);
    \TorusSphereApart{20.5}
\end{tikzpicture}}
  \caption{Illustration of Theorem \ref{tunnel counterexample}.}
    \label{figure counterexample}
\end{figure}
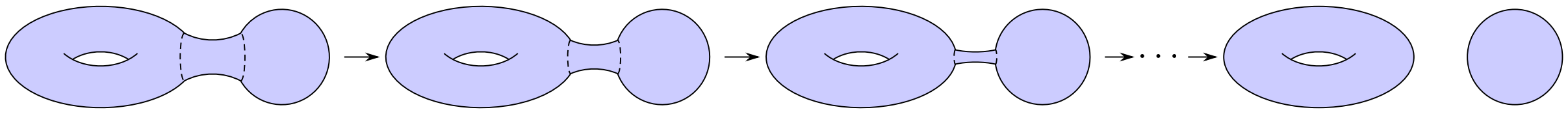

We apply this construction to disconnected Riemannian manifolds where one component is $S^2\times S^1$, which will be responsible to ensure the almost equality in the Bray-Brendle-Neves theorem, and the second component is $S^3$. The reason we choose to attach the tunnel with a $S^3$-factor is to ensure that the 2-spheres given by the slices of the tunnel that potentially could decrease $\mathcal{A}_{\pi_2}$ will not be admissible competitors for this infimum, since they will represent a trivial class in $\pi_2$. The resulting sequence remains connected, but using Theorem \ref{disjoint Lakzian-Sormani} we are able to show that the tunnels collapse in the intrinsic flat limit, producing a disconnected manifold and, thereby, obstructing the expected rigidity conclusion. Figure \ref{figure counterexample} illustrates the phenomenon that is happening in Theorem \ref{tunnel counterexample}.

We first prove Theorem \ref{disjoint Lakzian-Sormani}, which has its own interest and it is not related to scalar curvature.

\begin{proof}[Proof of Theorem \ref{disjoint Lakzian-Sormani}]
We prove the result directly from the definition of the intrinsic flat distance by constructing a common complete metric space and suitable integral currents. Fix a positive number $b$ such that $2b\geq\lambda$. Consider the disjoint union \[W\coloneqq M\sqcup X\] and define a metric $d_W$ on $W$ by preserving the original distances on each component and setting for $z\in M$ and $w\in X$,
\begin{equation}\label{bridge metric}
d_W(z,w)\coloneqq \inf_{x\in U} \left\{d_g(z,F(x))+b+d_X(x,w) \right\}.
\end{equation}
Also, put $d_W(w,z)\coloneqq d_W(z,w)$. Let $K \colon W\longrightarrow Z\coloneqq\ell^\infty(W) $ be the Kuratowski isometric embedding and define \[\varphi \coloneqq K\vert_M, \qquad \psi\coloneqq K\vert_X\] which, by definition of $d_W$, are isometric embeddings. And for every $y\in U$, choosing $x = y$ in the infimum defining \eqref{bridge metric}, we obtain the relation
\begin{equation}\label{phi psi b}
\Vert\varphi(F(y))-\psi(y)\Vert_Z = d_W(F(y),y) \leq b.
\end{equation}

Since $U=U_1\sqcup U_2$ is disconnected, we construct the currents separately on each component. For $i \in \{1,2\}$, equip $U_i$ with its intrinsic Riemannian distance associated with $\hat{h}_i$, denoted by $d_{\hat{U}_i}$, and set $Q_i\coloneqq\llbracket U_i\rrbracket \in\mathbf I_n(U_i).$ Define the linear interpolation $H_i\colon [0,1]\times U_i\longrightarrow Z$ by \[H_i(t,x)\coloneqq (1-t) \varphi(F(x)) + t\psi(x).\] We first establish uniform Lipschitz estimates for $H_i$. Since $F$ preserves the Riemannian metric tensor on $U_i$, it preserves the lengths of curves contained in $U_i$. Therefore, for every $x,y\in U_i$, \[\Vert \varphi(F(x))-\varphi(F(y))\Vert_Z =d_g(F(x),F(y)) \leq d_{\hat{U}_i}(x,y). \] Similarly, since $\psi$ is an isometric embedding and $d_X$ restricts to $d_{h_i}$ on $N_i$, we have \[\Vert \psi(x)-\psi(y)\Vert_Z =d_X(x,y) =d_{h_i}(x,y) \leq L_i d_{\hat{U}_i}(x,y)\] where the last inequality follows from the assumption $h_i(v,v) \leq L_i^2 \hat{h}_i(v,v)$ for all $v\in TU_i$.

Therefore, by the definition of $H_i$, \[\Vert H_i(t,x)-H_i(t,y)\Vert_Z \leq (1-t)\Vert \varphi(F(x)) -\varphi(F(y))\Vert_Z + t\Vert\psi(x)-\psi(y)\Vert_Z \leq L_i d_{\hat{U}_i}(x,y).\] Thus, for every $t \in[0,1]$, $\text{Lip}(H_i(t,\cdot))\leq L_i$. On the other hand, by \eqref{phi psi b},
\[ \Vert H_i(t,x)-H_i(s,x)\Vert_Z = \vert t-s \vert \Vert \psi(x)-\varphi(F(x))\Vert_Z \leq b\vert t-s\vert\] and hence $ \text{Lip}(H_i(\cdot,x))\leq b$ for every $x\in U_i.$ In particular, we proved that $H_i$ is Lipschitz on $[0,1]\times U_i$ since we endowed this space with the product distance. Also, the Lipschitz constants satisfy
\begin{equation}\label{H Lipschitz}
\text{Lip}(H_i(t,\cdot))\leq L_i, \qquad \text{Lip}(H_i(\cdot,x))\leq b
\end{equation}
for every $t\in[0,1]$ and $x\in U_i$.

Define the $(n+1)$-dimensional integral current \[B_i\coloneqq (H_i)_\#([0,1]\times Q_i) \in\mathbf I_{n+1}(Z) \] and the $n$-dimensional integral current which is part of its boundary \[S_i\coloneqq (H_i)_\#([0,1]\times\partial Q_i) \in\mathbf I_n(Z).\] Set $B\coloneqq B_1+B_2$ and  $S\coloneqq S_1+S_2$. We now estimate the masses of the currents $B$ and $S$. For this, we follow the computations given in Wenger \cite[Theorem 2.9]{Wenger2005} used in the proof of Sormani-Wenger \cite[Lemma 5.4]{SormaniWenger2011}.

Let $(f,\pi_1,\ldots,\pi_{n+1}) \in\mathcal D^{n+1}(Z)$ be a test tuple. By the definition of the pushforward current, \[B_i(f,\pi_1,\ldots,\pi_{n+1}) = ([0,1]\times Q_i) (f\circ H_i,\pi_1\circ H_i,\ldots, \pi_{n+1}\circ H_i).\] For each $t\in[0,1]$, write $H_{i,t}(\cdot)=H_i(t,\cdot)$ and $\pi_{k,t}(\cdot) = \pi_k(H_i(t,\cdot)).$ Using the definition of the product current given in \cite[Definition 2.8]{Wenger2005}, we obtain
\begin{equation}\label{product current}
B_i(f,\pi_1,\ldots,\pi_{n+1}) = \sum_{k=1}^{n+1}(-1)^{k+1} \int_0^1
Q_i\left((f\circ H_{i,t}) \frac{\partial\pi_{k,t}}{\partial t}, \pi_{1,t},\ldots, \widehat{\pi_{k,t}},\ldots, \pi_{n+1,t} \right)\ dt,
\end{equation}
where the hat denotes omission.

Recall that the mass inequality for the $n$-dimensional current $Q_i$ gives for $(g, \sigma_1, \dots,\sigma_n) \in \mathcal D^n(U_i)$, \[\vert Q_i(g,\sigma_1,\ldots,\sigma_n)\vert \leq \prod_{\ell=1}^n\text{Lip}(\sigma_\ell) \int_{U_i}\vert g(x)\vert\ d\Vert Q_i\Vert(x). \] We apply this inequality to each term in \eqref{product current}. First, by \eqref{H Lipschitz}, \[\text{Lip}(\pi_{k,t}) = \text{Lip}(\pi_k\circ H_{i,t}) \leq \text{Lip}(\pi_k) \ \text{Lip}(H_{i,t}) \leq L_i \text{Lip}(\pi_k).\] Furthermore, for every $x\in U_i$ and $s,t\in[0,1]$, the estimate \eqref{H Lipschitz} also gives \[\vert \pi_{k,t}(x)-\pi_{k,s}(x)\vert  \leq \text{Lip}(\pi_k) \Vert H_i(t,x)-H_i(s,x)\Vert_Z \leq b \ \text{Lip}(\pi_k)\vert t-s\vert.\] Thus, the function $t\mapsto\pi_{k,t}(x)$ is Lipschitz. Therefore, it is differentiable for almost every $t$ and \[ \left\vert \frac{\partial\pi_{k,t}}{\partial t}(x) \right\vert \leq b\ \text{Lip}(\pi_k)\] holds for almost every $(t,x)$. Combining these inequalities with the mass inequality for $Q_i$, we obtain
\begin{align*} &\Bigg\vert Q_i\left(
(f\circ H_{i,t}) \frac{\partial\pi_{k,t}}{\partial t}, \pi_{1,t},\ldots, \widehat{\pi_{k,t}},\ldots, \pi_{n+1,t} \right) \Bigg\vert \\
&\qquad\leq \prod_{\ell\neq k} \text{Lip}(\pi_{\ell,t}) \int_{U_i} \vert f(H_i(t,x))\vert  \left\vert  \frac{\partial\pi_{k,t}}{\partial t}(x) \right\vert  d\Vert Q_i\Vert (x)\\
&\qquad\leq b L_i^n \prod_{\ell=1}^{n+1} \text{Lip}(\pi_\ell) \int_{U_i} \vert f(H_i(t,x))\vert \ d\Vert Q_i\Vert(x).
\end{align*}
Replacing this estimate into \eqref{product current} and using the triangle inequality, we conclude that \begin{equation}\label{product current 2}
\vert B_i(f,\pi_1,\ldots,\pi_{n+1})\vert \leq (n+1)b L_i^n\prod_{\ell=1}^{n+1}\text{Lip}(\pi_\ell) \int_0^1\int_{U_i} \vert f(H_i(t,x))\vert \ d\Vert Q_i\Vert(x)\ dt.
\end{equation}

Consider now the finite Borel measure on $Z$ defined by \[\mu_i\coloneqq (n+1)bL_i^n (H_i)_\# \bigl(\mathcal L^1\llcorner[0,1] \otimes \Vert Q_i\Vert \bigr).\] By the definition of the pushforward of a measure, \eqref{product current 2} can be rewritten as \[\vert B_i(f,\pi_1,\ldots,\pi_{n+1})\vert  \leq \prod_{\ell=1}^{n+1} \text{Lip}(\pi_\ell) \int_Z\vert f\vert \ d\mu_i.\] By the definition of the mass measure, we conclude that $\Vert B_i\Vert \leq\mu_i.$ In particular,
\begin{align*}
\mathbf M(B_i) = \Vert B_i\Vert (Z) \leq \mu_i(Z) &= (n+1)b L_i^n \left(\mathcal L^1 \llcorner[0,1] \otimes\Vert Q_i\Vert \right)([0,1]\times U_i)\\
&=(n+1)bL_i^n\ \Vert Q_i\Vert (U_i) =(n+1)b L_i^n\ \mathbf M(Q_i)\\
&=(n+1)bL_i^n \vol_{\hat{h}_i}(U_i)
\end{align*}
where the last equality holds since $\mathbf M(Q_i)= \vol_{\hat{h}_i}(U_i)$ given that $Q_i=\llbracket U_i\rrbracket$. 

Now recall that $S_i=(H_i)_\#([0,1]\times\partial Q_i).$ Since $\partial Q_i$ is an $(n-1)$-dimensional integral current, a similar calculation yields \[\mathbf M(S_i) \leq nb L_i^{n-1}\ \mathbf M(\partial Q_i) = nb L_i^{n-1}\ \text{Area}_{\hat{h}_i}(\partial U_i).\] Finally, since $B=B_1+B_2$ and $S=S_1+S_2$, the triangle inequality for the mass gives \[\mathbf M(B) \leq\sum_{i=1}^2\mathbf M(B_i) \leq(n+1)b \sum_{i=1}^2 L_i^n \vol_{\hat{h}_i}(U_i)\] and \[\mathbf M(S) \leq\sum_{i=1}^2\mathbf M(S_i) \leq n b \sum_{i=1}^2 L_i^{n-1} \text{Area}_{\hat{h}_i}(\partial U_i).\]

Recall the boundary formula for the product current \[\partial([0,1]\times Q_i) = [1]\times Q_i -[0]\times Q_i -[0,1]\times\partial Q_i.\] Since pushforward commutes with boundary, $H_i(0,\cdot) =\varphi (F(\cdot))$ and $H_i(1,\cdot) = \psi(\cdot)$, we obtain \[ \partial B_i = (\psi\vert_{U_i})_\#Q_i - (\varphi\circ F\vert_{U_i})_\#Q_i -S_i.\] Summing over $i$, we conclude that
\begin{equation}\label{boundary B}
\partial B = \sum_{i=1}^2(\psi\vert_{U_i})_\#Q_i -\sum_{i=1}^2(\varphi\circ F\vert_{U_i})_\#Q_i -S.
\end{equation}

Since the boundaries of the regions $U_i$ have zero $n$-dimensional volume, the currents $T_M$ and $T_X$ given by integration may be decomposed as
\[T_M = \sum_{i=1}^2 (F\vert_{U_i})_\#Q_i + T_M\llcorner(M\setminus F(U)) \] and \[ T_X = T_X\llcorner U + T _X\llcorner(X\setminus U).\]
Defining the $n$-dimensional integral current \[ A\coloneqq -S+\varphi_\#\big(T_M\llcorner (M\setminus F(U))\big)-\psi_\#(T_X\llcorner (X\setminus U))\] and combining the preceding decompositions with \eqref{boundary B}, we obtain \[\varphi_\#T_M-\psi_\#T_X =\partial(-B)+A.\] Furthermore, since $\varphi$ and $\psi$ are isometric embeddings, the triangle inequality for mass yields \[ \mathbf M(A) \leq \mathbf M(S) + \vol _g(M\setminus F(U))+\mathcal H^n_{d_X}(X\setminus U).\] By the definition of the intrinsic flat distance,
\begin{align*}
d_{\mathcal F} \bigl((M,d_g,T_M),(X, d_X,& \ T_X)\bigr)\leq \mathbf M(B)+\mathbf M(A) \leq (n+1)b \sum_{i=1}^2 L_i^{n}\vol_{\hat{h}_i}(U_i)\\
& +nb\sum_{i=1}^2L_i^{n-1} \text{Area}_{\hat{h}_i}(\partial U_i) + \vol_g(M\setminus F(U)) +\mathcal H^n_{d_X}(X\setminus U).
\end{align*}
Since this estimate holds for every positive number $b$ satisfying $2b\geq\lambda$, we may let $b\downarrow\lambda/2$, obtaining
\begin{align*}
d_{\mathcal F} \bigl((M,d_g,T_M),(X,d_X,T_X)\bigr) \leq \frac{n+1}{2} \lambda \sum_{i=1}^2 \ & L_i^{n}\vol_{\hat{h}_i}(U_i)+
\frac{n}{2}\lambda \sum_{i=1}^2 L_i^{n-1}\text{Area}_{\hat{h}_i}(\partial U_i)\\
&+\vol_g(M\setminus F(U)) + \mathcal H^n_{d_X}(X\setminus U).
\end{align*}
\end{proof}

\begin{rmk}
    Theorem \ref{disjoint Lakzian-Sormani} can be used to prove Proposition 5.1 in \cite{Sweeney2025}.
\end{rmk}

\begin{rmk}
    Theorem \ref{disjoint Lakzian-Sormani} can be adapted to the case where a tunnel is attached between two distinct points $p,q$ of the same connected manifold $(N,h)$, creating a handle. In this case, we let $Y$ be the length space obtained by attaching a segment of length $d>0$ to $N$ at $p$ and $q$. We then define $X\coloneqq  N$,  $d_X \coloneqq d_Y\vert_{N\times N}.$ Unlike the case of two disjoint manifolds, $d_X$ may be strictly smaller than $d_h$. Nevertheless, the two metrics agree locally, so the mass of the integral current induced by $N$ is still equal to $\vol_h(N)$. The proof then follows the same argument, taking $U$ to be the complement of the two balls removed from $N$.
\end{rmk}

We now prove Theorem \ref{tunnel counterexample}. We divide it into two parts: in Part I we verify that the choice of manifolds explained in the beginning of this section satisfies the almost equality condition of Theorem \ref{BBN theorem}. In Part II we apply Theorem \ref{disjoint Lakzian-Sormani} to guarantee the volume preserving intrinsic flat convergence.

\begin{proof}[Proof of Theorem \ref{tunnel counterexample}]

\textbf{Part I.} 
Pick the sequence of radii $r_j\coloneqq \sqrt{\frac{j}{1+j}}$ and for each $j \in \mathbb{N}$ consider the product manifold \[X_{1,j} \coloneqq (X_1,h_{1,j}) = (S^2 \times S^1_\ell, \ r_j^2g_{S^2}+dt^2) \qquad \text{and} \qquad X_2 \coloneqq (\mathbb S^3,h_2) = (\mathbb S^3, g_{S^3})\] where $g_{S^2}$, $g_{S^3}$ are the round metrics and $\ell>0$ is the $S^1$-length. Then the scalar curvature of $X_{1,j}$ is constant $R_{h_{1,j}} \equiv 2+\frac{2}{j}.$ Notice also that each $X_{1,j}$ has $\mathcal{A}_{\pi_2}(X_{1,j},h_{1,j}) = 4\pi r_j^2$, since this infimum is realized by any slice $S^2\times \{t\}$.

Choose points $q \in X_2$, $t_* \in S^1_\ell$ and $p \in X_1$ which does not lie in the slice $\Sigma_j \coloneq S^2\times \{t_*\}$. Fix $d>0$ and for the Riemannian manifolds $X_{1,j}$ and $X_2$, consider $\delta_j >0$  so small that Proposition \ref{tunnel} applies with $\delta_j < \frac{1}{j}$, the ball of radius $4\sqrt\delta_j$ about the point $p$ misses $\Sigma_j$ and for every $j$ this ball lies in one fixed normal coordinate neighborhood $\mathcal U \coloneqq B_R^{h_1}(p)$ with respect to the limit metric $h_1 \coloneqq g_{S^2}+dt^2$ where $0< R < \text{inj}_{h_1}(p)$. Denote $(M_j,g_j)$ the resulting manifold. Topologically $M_j$ is the connected sum $X_1 \# S^3$ and so it is diffeomorphic to $S^2\times S^1$. Also the scalar curvature of $M_j$ satisfies \[R_{g_j} > 2+\dfrac{2}{j}-\dfrac{1}{j} = 2+\dfrac{1}{j}.\] Since outside the surgery region the metric is still $h_{1,j}$, there are points where the scalar curvature of $M_j$ is exactly $2+\frac{2}{j}$. Hence the inequalities \[2+\dfrac{1}{j} \leq \inf_{M_j} R_{g_j} \leq 2+\dfrac{2}{j}\] hold, and so $\inf_{M_j} R_{g_j} \xrightarrow{\ j \rightarrow +\infty \ } 2.$ We claim that \[\mathcal A_{\pi_2}(M_j,g_j) \xrightarrow{\ j \rightarrow +\infty \ } 4\pi.\] 

Indeed, because the tunnel is attached outside $\Sigma_j$ and the metric remains unchanged there, then $\mathcal{A}_{\pi_2}(M_j,g_j) \leq \text{Area}_{h_{1,j}}(\Sigma_j) = 4\pi r_j^2.$ Regarding the lower bound, we denote $S^2(r_j)$ the round $2$-sphere with radius $r_j$ and in $X_{1,j}$ consider the $2$-form $\omega_j$ given by \[\omega_j = \pi^* dA_{S^2(r_j)}\] where $dA_{S^2(r_j)}$ is the area element of $S^2(r_j)$ and $\pi$ is the projection of $X_{1,j}$ into $S^2(r_j)$. By construction \[\int_{\Sigma_j}\omega_j = 4\pi r_j^2.\] We claim that its comass is one. In fact, the comass of $\omega_j$ is defined as \[\Vert \omega_j \Vert_{*,h_{1,j}} = \sup\{ \vert \omega_j(e_1,e_2)\vert \colon \vert e_1 \vert = \vert e_2 \vert = 1 \quad \text{and} \quad e_1 \perp e_2\}.\] Fix $(x,t) \in X_{1,j}$ and take two orthonormal vectors $v_1,v_2 \in T_{(x,t)}X_{1,j}$. Then \[\omega_j(v_1,v_2) =dA_{S^2(r_j)}(d\pi(v_1), d\pi(v_2))\] and since $\pi$ is a Riemannian submersion, $\vert d\pi(v_i)\vert \leq 1$ and so $\vert \omega_j(v_1,v_2)\vert \leq 1.$ Hence $\Vert \omega_j \Vert_{*,h_{1,j}} \leq 1$. To prove the equality, consider $\hat{e}_1, \hat{e}_2$ positively oriented orthonormal basis of $T_xS^2(r_j)$. We can see them as vectors in $T_{(x,t)}X_{1,j}$ by setting $e_1 = (\hat{e}_1,0)$ and $e_2 = (\hat{e}_2,0)$. In the product metric $h_{1,j}$ they are orthonormal and so \[\omega_j(e_1,e_2) = dA_{S^2(r_j)}(d\pi (e_1), d\pi(e_2)) = dA_{S^2(r_j)}(\hat{e}_1, \hat{e}_2) = 1.\] 

Now we want to define a closed 2-form on the new manifold $M_j$ which equals $\omega_j$ away from a neighborhood of the surgery in $X_{1,j}$ and which vanishes on the tunnel and in the $S^3$ part. By the definition of $\mathcal U$, for each $j$ the proof of Poincaré Lemma in \cite[Section 4.18]{Warner1983} guarantees the existence of a $1$-form $\alpha_j$ there such that $\omega_j = d\alpha_j$ and \[\Vert \alpha_j(x)\Vert_{*,h_{1,j}} \leq C \ d_{h_{1,j}}(p,x)\] for all $x$ in the $h_{1,j}$-ball  $B_{4\sqrt{\delta_j}}(p)$. The constant $C>0$ can be chosen independently of $j$, since $h_{1,j}\to h_1$ smoothly and the coordinate coefficients of $\omega_j$ are uniformly bounded.

For $\rho_j \coloneqq 3\sqrt{\delta_j}$, $\rho_j > 3\delta_j$ and $\frac{\rho_j}{\delta_j} \rightarrow +\infty$. Choose a radial cut-off function $f_j=f_j(r)$ on $M_j$ such that 
\[f_j = \begin{cases}
        0 & \text{for } r \leq 3\delta_j\\
        1 & \text{for } r \geq \rho_j
        \end{cases}\]
and between these two values $f_j$ satisfies \[\vert df_j \vert \leq\dfrac{K}{r \log(\frac{\rho_j}{3\delta_j})}\] for some constant $K>0$. We can then define a smooth 2-form $\Omega_j$ in $M_j$ by \[\Omega_j \coloneqq \begin{cases}
    \omega_j, &r\geq \rho_j \\
    d(f_j\alpha_j), & 3\delta_j < r < \rho_j\\
    0 & \text{otherwise.}
\end{cases}\]
By construction, $d\Omega_j = 0$. Also, in the transition annulus we have \[\Omega_j =d(f_j\alpha_j) =f_j \omega_j + df_j \wedge \alpha_j\] and so  \[\Vert \Omega_j\Vert_{*,g_j} \leq 1+\dfrac{K}{\log(\frac{\rho_j}{3\delta_j})}\] in this region. Since $\rho_j$ was chosen in a way that $\frac{\rho_j}{\delta_j} \rightarrow +\infty$, then there exists a sequence of positive numbers $\eta_j \rightarrow 0$ such that
\begin{equation}\label{estimate Omega}
    \Vert \Omega_j\Vert_{*,g_j} \leq 1+\eta_j.    
\end{equation}
Actually the estimate \eqref{estimate Omega} holds in all regions, since $\Vert \omega_j\Vert_{*,h_{1,j}} = 1$. 

We now explain how $\Vert \Omega_j\Vert_{*,g_j}$ is related to $\mathcal{A}_{\pi_2}(M_j, g_j)$. Take $f\colon S^2 \rightarrow M_j$ a smooth map such that its homotopy class $[f]\in \pi_2(M_j)$ is nonzero. As $M_j \simeq (S^2\times S^1) \# S^3$, $[f]\neq 0$ and $\pi_2(S^1) =0$, then $f$ is homotopic to a map whose image is contained in $\Sigma_j$. Since $\Sigma_j$ is far away from the ball where the surgery occurred, then $\Omega_j\vert_{\Sigma_j} = \omega_j$ and so \[\int_{S^2}f^*\Omega_j = 4\pi mr_j^2\] for some nonzero integer $m$. Therefore, 
\begin{equation}\label{modulo da integral fOmega}
    \left\vert \int_{S^2}f^*\Omega_j \right\vert = 4\pi \vert m \vert r_j^2 \geq 4\pi r_j^2.
\end{equation}
On the other hand, let $x \in S^2$ and take $e_1,e_2 $ a positively oriented orthonormal basis of $T_xS^2$. Then \[(f^*\Omega_j)_x(e_1,e_2) = (\Omega_j)_{f(x)}(df_x(e_1), df_x(e_2))\] and by the definition of comass, \[\left\vert (\Omega_j)_{f(x)}(df_x(e_1),df_x(e_2))\right\vert \leq \Vert \Omega_j\Vert_{*,g_j} \vert df_x(e_1) \wedge df_x(e_2)\vert_{g_j} = \Vert \Omega_j\Vert_{*,g_j} J_2f(x) \] where $J_2f$ is the Jacobian of $f$. Therefore, by the area formula, \[\left\vert \int_{S^2}f^*\Omega_j\right\vert \leq \int_{S^2}\left\vert f^*\Omega_j\right\vert \leq \Vert \Omega_j\Vert_{*,g_j} \int_{S^2}J_2f \ dA =  \Vert \Omega_j\Vert_{*,g_j} \ \text{Area}_{f^*g_j}(S^2).\] Putting this estimate on (\ref{modulo da integral fOmega}), \[4\pi r_j^2 \leq \Vert \Omega_j\Vert_{*,g_j} \ \text{Area}_{f^*g_j}(S^2) \leq (1+\eta_j)\ \text{Area}_{f^*g_j}(S^2)\] and then \[\text{Area}_{f^*g_j}(S^2) \geq \dfrac{4\pi r_j^2}{1+\eta_j}.\] Taking the infimum over all nonzero $[f] \in \pi_2(M_j)$, \[\dfrac{4\pi r_j^2}{1+\eta_j} \leq \mathcal{A}_{\pi_2}(M_j,g_j) \leq 4\pi r_j^2.\] Since $r_j \rightarrow1$ and $\eta_j \rightarrow 0$,  then $\mathcal{A}_{\pi_2}(M_j, g_j) \rightarrow 4\pi.$ Therefore the manifolds $(M_j, g_j)$ satisfy \[\mathcal{A}_{\pi_2}(M_j, g_j)\inf_{M_j} R_{g_j} \xrightarrow{\ j \rightarrow +\infty \ } 8\pi.\] 

\noindent\textbf{Part II.}
We now prove the volume preserving intrinsic flat convergence. Recall that \[(X_{1,j},h_{1,j}) = ( S^2\times S^1_\ell,\ r_j^2g_{S^2}+dt^2) \qquad \text{and} \qquad (X_2,h_2)=(\mathbb S^3,g_{S^3}).\] Also recall that the Riemannian manifold $(M_j,g_j)$ was obtained by removing the geodesic balls $B_{2\delta_j}^{h_{1,j}}(p)$ and $B_{2\delta_j}^{h_2}(q)$ and joining a tunnel of prescribed central length $d>0$. By construction, the tunnel can be decomposed (see Figure \ref{figure tunnel} for guidance) as $T_j=A_{1,j}\cup C_j\cup A_{2,j}$ where $C_j= \left([0,d]\times S^2, dt^2+c_j^2g_{ S^2} \right)$, with $c_j \rightarrow 0$, and the half tunnels $A_{1,j}$, $A_{2,j}$ satisfy (see Section $3.6$ in \cite{Sweeney2025}) the estimates \[\text{diam}(A_{i,j})\leq C\delta_j \qquad \text{and} \qquad \vol(A_{i,j})\leq C\delta_j^3\] for some constant $C>0$. The constant $C$ in principle may depend on the underlying metric $h_{1,j}$, but since $h_{1,j} \rightarrow h_1 \coloneqq g_{S^2}+dt^2$ smoothly and the background manifolds have uniformly controlled geometry, this constant can be chosen uniformly in $j$.
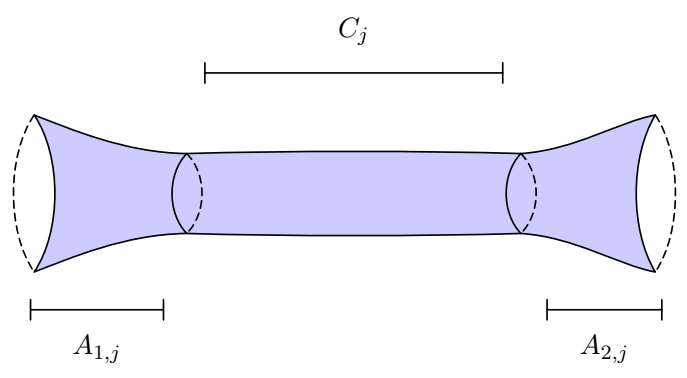
\begin{figure}[H]
\centering
\resizebox{0.28\linewidth}{!}{%
\begin{tikzpicture}[
    line width=0.65pt,
    line cap=round,
    line join=round
]

\colorlet{myfill}{blue!20}

\filldraw[fill=myfill, draw=black]
  (-4.4,1.08)
  .. controls (-3.90,0.89) and (-3.15,0.56)
  .. (-2.30,0.55)
  .. controls (-0.95,0.59) and (0.95,0.59)
  .. (2.30,0.55)
  .. controls (3.05,0.60) and (3.65,0.99)
  .. (4.15,1.08)
  .. controls (3.80,0.50) and (3.80,-0.50)
  .. (4.15,-1.08)
  .. controls (3.65,-0.99) and (3.05,-0.60)
  .. (2.30,-0.55)
  .. controls (0.95,-0.59) and (-0.95,-0.59)
  .. (-2.30,-0.55)
  .. controls (-3.15,-0.56) and (-3.90,-0.89)
  .. (-4.4,-1.08)
  .. controls (-4.02,-0.52) and (-4.02,0.52)
  .. (-4.4,1.08)
  -- cycle;

\draw[densely dashed]
  (-4.4,1.08)
  .. controls (-4.78,0.53) and (-4.78,-0.53)
  .. (-4.4,-1.08);

\draw[densely dashed]
  (4.15,1.08)
  .. controls (4.52,0.53) and (4.52,-0.53)
  .. (4.15,-1.08);

\draw
  (-2.30,0.55)
  .. controls (-2.57,0.27) and (-2.57,-0.27)
  .. (-2.30,-0.55);

\draw[densely dashed]
  (-2.30,0.55)
  .. controls (-2.02,0.27) and (-2.02,-0.27)
  .. (-2.30,-0.55);

\draw
  (2.30,0.55)
  .. controls (2.03,0.27) and (2.03,-0.27)
  .. (2.30,-0.55);

\draw[densely dashed]
  (2.30,0.55)
  .. controls (2.58,0.27) and (2.58,-0.27)
  .. (2.30,-0.55);

\draw (-2.05,1.66) -- (2.05,1.66);
\draw (-2.05,1.52) -- (-2.05,1.80);
\draw (2.05,1.52) -- (2.05,1.80);
\node[above=7pt] at (0,1.66) {$C_j$};

\draw (-4.45,-1.60) -- (-2.62,-1.60);
\draw (-4.45,-1.74) -- (-4.45,-1.46);
\draw (-2.62,-1.74) -- (-2.62,-1.46);
\node[below=6pt] at (-3.535,-1.60)
  {$A_{1,j}$};

\draw (2.66,-1.60) -- (4.24,-1.60);
\draw (2.66,-1.74) -- (2.66,-1.46);
\draw (4.24,-1.74) -- (4.24,-1.46);
\node[below=6pt] at (3.45,-1.60)
  {$A_{2,j}$};

\end{tikzpicture}}
  \caption{Illustration of the tunnel.}
    \label{figure tunnel}
\end{figure}

Let $(Y,d_Y)$ be the length space obtained by joining $(X_1,h_1)$ and $(X_2,h_2)$ with a segment of length $d$ attached at $p$ and $q$, respectively (see Figure \ref{figure Y} for guidance), and equipped with the induced length metric. Define \[X \coloneqq X_1\sqcup X_2, \qquad d_X \coloneqq d_Y|_{X\times X}\] and let $\mathcal T\coloneqq \llbracket X_1\rrbracket + \llbracket X_2\rrbracket $ be the corresponding integral current with orientations compatible with those of $M_j$. 

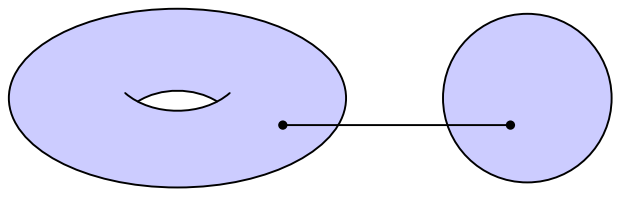
\begin{figure}[H]
\centering
\resizebox{0.25\linewidth}{!}{%
\begin{tikzpicture}[
    line width=0.65pt,
    line cap=round,
    line join=round
]

\colorlet{myfill}{blue!20}

\newcommand{\TorusHolePathh}{%
  (-0.62,0.06)
  .. controls (-0.3,-0.22) and (0.3,-0.22) ..
  (0.62,0.06)
  -- (0.47,-0.04)
  .. controls (0.2,0.13) and (-0.2,0.13) ..
  (-0.47,-0.04)
  -- cycle
}

\newcommand{\TorusHolee}{%
  \draw (-0.62,0.06)
    .. controls (-0.3,-0.22) and (0.3,-0.22) ..
    (0.62,0.06);

  \draw (-0.47,-0.04)
    .. controls (-0.2,0.13) and (0.2,0.13) ..
    (0.47,-0.04);
}

\begin{scope}[xshift=-2.3cm]

  \path[fill=myfill, even odd rule]
    (0,0) ellipse[x radius=2, y radius=1.06]
    \TorusHolePathh;

  \draw (0,0)
    ellipse[x radius=2, y radius=1.06];

  \TorusHolee

\end{scope}

\filldraw[fill=myfill, draw=black]
  (1.85,0) circle[radius=1];

\draw (-1.05,-0.32) -- (1.65,-0.32);

\fill (-1.05,-0.32) circle (0.055);
\fill (1.65,-0.32) circle (0.055);

\end{tikzpicture}}
  \caption{Illustration of the length space $Y = X_1 \cup_p [0,d] \cup_qX_2$.}
    \label{figure Y}
\end{figure}

We claim that \[(M_j,d_{g_j},\llbracket M_j\rrbracket) \xrightarrow{\ \mathcal F\ } (X,d_X,\mathcal T).\] For that, we verify the hypotheses of Theorem \ref{disjoint Lakzian-Sormani}. Consider the regions \[U_{1,j} \coloneqq X_1\setminus B_{2\delta_j}^{h_{1,j}}(p), \qquad U_{2,j} \coloneqq X_2\setminus B_{2\delta_j}^{h_2}(q)\] and set $ U_j \coloneqq U_{1,j}\sqcup U_{2,j}\subset X$ (see Figure \ref{figure U} for guidance). Let $F_j \colon U_j \longrightarrow M_j$ be the natural identification of the unchanged regions with their images in $M_j$.

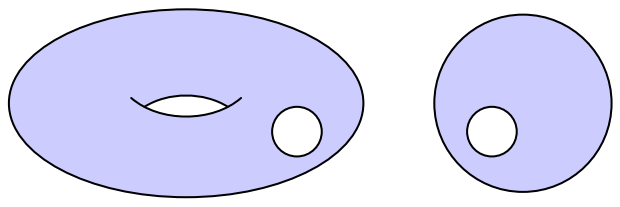
\begin{figure}[H]
\centering
\resizebox{0.25\linewidth}{!}{%
\begin{tikzpicture}[
    line width=0.65pt,
    line cap=round,
    line join=round
]

\colorlet{myfill}{blue!20}

\newcommand{\TorusHolePathhh}{%
  (-0.62,0.06)
  .. controls (-0.3,-0.22) and (0.3,-0.22) ..
  (0.62,0.06)
  -- (0.47,-0.04)
  .. controls (0.2,0.13) and (-0.2,0.13) ..
  (-0.47,-0.04)
  -- cycle
}

\newcommand{\TorusHoleee}{%
  \draw (-0.62,0.06)
    .. controls (-0.3,-0.22) and (0.3,-0.22) ..
    (0.62,0.06);

  \draw (-0.47,-0.04)
    .. controls (-0.2,0.13) and (0.2,0.13) ..
    (0.47,-0.04);
}

\coordinate (pT) at (1.25,-0.32);
\coordinate (pS) at (3.45,-0.32);

\def\r{0.28}


\path[fill=myfill, even odd rule]
    (0,0) ellipse[x radius=2, y radius=1.06]
    \TorusHolePathhh
    (pT) circle[radius=\r];

\draw (0,0)
    ellipse[x radius=2, y radius=1.06];

\TorusHoleee

\draw
    (pT) circle[radius=\r];


\path[fill=myfill, even odd rule]
    (3.8,0) circle[radius=1]
    (pS) circle[radius=\r];

\draw (3.8,0)
    circle[radius=1];

\draw
    (pS) circle[radius=\r];

\end{tikzpicture}}
  \caption{Illustration of the region $U_j$.}
    \label{figure U}
\end{figure}

Since $h_{1,j}\rightarrow h_1$ smoothly and the background manifolds have uniformly controlled geometry, we have the estimates
\begin{equation}\label{hausdorff measure of X U}
    \mathcal H^3_{d_X}(X\setminus U_j) \leq C\delta_j^3 \qquad \text{and} \qquad \mathcal H^2_{d_X}(\partial U_j) \leq C\delta_j^2    
\end{equation}
for a constant $C>0$ that does not depend on $j$. Here $\partial U_j$ is the disjoint union of the two geodesic spheres and the Hausdorff measures are computed with respect to $d_X$. Furthermore, $M_j\setminus F_j(U_j)=T_j$ up to sets of zero $3$-dimensional measure. By the assumptions on the tunnel, \[\vol(T_j) \leq \vol(A_{1,j}) + \vol(C_j) +\vol(A_{2,j})\] which implies that 
\begin{equation}\label{volume of M U}  
    \vol_{g_j} (M_j\setminus F_j(U_j)) \longrightarrow 0.
\end{equation}

We claim that there exists a sequence $\varepsilon_j\rightarrow 0$ such that
\begin{equation}\label{distance distorcion}
\lambda_j \coloneqq \sup_{x,y\in U_j} \vert d_{g_j}(F_j(x),F_j(y))-d_X(x,y) \vert \leq \varepsilon_j.
\end{equation}

In fact, write $d_{1,j}$, $d_1$ and $d_2$ for the Riemannian distances associated with $h_{1,j}$, $h_1$ and $h_2$, respectively. First suppose that $x,y\in U_{1,j}$. Since the intrinsic diameter of $\partial B_{2\delta_j}^{h_{1,j}}(p)$ is bounded by $C\delta_j$, a minimizing path in $X_1$ can be modified to avoid the removed ball at an additional length of at most $C\delta_j$. Conversely, any path in $M_j$ joining $F_j(x)$ to $F_j(y)$ that enters the tunnel can be modified to remain in $X_1$ by replacing the portion between its first and last intersections with the attachment sphere by a path along that sphere.

Therefore, \[ \vert d_{g_j}(F_j(x),F_j(y))-d_{1,j}(x,y) \vert \leq C\delta_j.\] Since $h_{1,j}\rightarrow h_1$ smoothly, we have \[a_j\coloneqq  \sup_{x,y\in X_1} \vert d_{1,j}(x,y)-d_1(x,y)\vert \longrightarrow 0\] and with this it follows that
\begin{equation}\label{dist 1}
\vert d_{g_j}(F_j(x),F_j(y))-d_1(x,y) \vert \leq C\delta_j+ a_j.
\end{equation}

The same argument on the second component gives
\begin{equation}\label{dist 2}
\vert d_{g_j}(F_j(x),F_j(y))-d_2(x,y) \vert \leq C\delta_j
\end{equation}
for every $x,y\in U_{2,j}$.

Now let $x\in U_{1,j}$ and $y\in U_{2,j}$. Every path joining $F_j(x)$ to $F_j(y)$ must pass through the central cylinder $C_j$. Hence the part of the path inside $C_j$ has length at least $d$. Thus \[d_{g_j}(F_j(x),F_j(y)) \geq d_{1,j}(x,p)+d +d_2(q,y)-4\delta_j.\] For the opposite inequality, we construct the following path: start at $x$, follow the minimizing geodesic toward $p$ until it reaches $\partial B_{2\delta_j}(p)$. Then traverse the first half tunnel $A_{1,j}$, the central cylinder $C_j$ along an axial curve and $A_{2,j}$. Finally follow a minimizing geodesic toward $y$. Since the diameter of $A_{i,j}$ for $i\in \{1,2\}$ is at most $C\delta_j$, we have: \[d_{g_j}(F_j(x),F_j(y)) \leq d_{1,j}(x,p)+ d +d_2(q,y)+C\delta_j.\] Therefore, \[\Big\vert d_{g_j}(F_j(x),F_j(y)) -d_{1,j}(x,p)-d-d_2(q,y) \Big\vert \leq C\delta_j.\] By the definition of $d_X$ we have $d_X(x,y)=d_1(x,p)+d+d_2(q,y)$, thus we conclude that
\begin{equation}\label{dist 3}
\vert d_{g_j}(F_j(x),F_j(y))-d_X(x,y)\vert \leq C\delta_j+a_j.
\end{equation}
Combining \eqref{dist 1}, \eqref{dist 2} and \eqref{dist 3}, we obtain \eqref{distance distorcion} with $\varepsilon_j=C\delta_j+a_j\longrightarrow 0$.

By Theorem \ref{disjoint Lakzian-Sormani}, choosing $L_{1,j} = r_j^{-1}$ and $L_{2,j} = 1$, for each $j$ it is possible to estimate the intrinsic flat distance as
\begin{align*}
d_{\mathcal F}
\bigl((M_j,d_{g_j},& \ \mathcal T_j),(X,d_X,\mathcal T) \bigr) \leq 2 \lambda_j\bigl(r_j^{-3}
\vol_{h_{1,j}}(U_{1,j})
+\vol_{h_2}(U_{2,j})
\bigr) \\
&+ \frac{3}{2}\lambda_j \bigl( r_j^{-2}
\text{Area}_{h_{1,j}}(\partial U_{1,j})
+\text{Area}_{h_2}(\partial U_{2,j})
\bigr)+ \vol_{g_j}(M_j \setminus F_j(U_j))
+\mathcal H^3_{d_X}(X\setminus U_j).
\end{align*}

By \eqref{hausdorff measure of X U}, \eqref{volume of M U}, \eqref{distance distorcion}, the fact that $\delta_j \rightarrow 0$ and $r_j \rightarrow 1$, we have that $(M_j,d_{g_j},\mathcal T_j) \xrightarrow{\ \mathcal F \ } (X,d_X,\mathcal T)$ as desired.

We claim that the volumes also converge. In fact, by the volume decomposition \[\vol_{g_j}(M_j) = \vol_{h_{1,j}}(X_{1,j}) -\vol_{h_{1,j}}(B_{2\delta_j}(p)) +\vol_{h_2}(X_2) -\vol_{h_2}(B_{2\delta_j}(q)) +\vol_{g_j}(T_j),\] since $r_j\rightarrow 1$, $\delta_j \rightarrow 0$ and $\vol(T_j)\rightarrow 0$, we obtain \[\lim_{j\to\infty}\vol_{g_j}(M_j) = 4\pi \ell +2\pi^2 = \vol_{h_1}(X_1) + \vol_{h_2}(X_2) =\mathbf M( \mathcal T).\] Hence the convergence is also volume preserving.
\end{proof}  

A related scalar curvature compactness conjecture deals with a smaller number related to a $3$-dimensional Riemannian Manifold $(M,g)$ given by \[\text{Min A}(M,g) \coloneqq \inf\{\text{Area}(\Sigma)  \colon  \Sigma \text{ is a closed minimal surface in } M^3\}.\] Observe that for Riemannian manifolds with nontrivial second homotopy group we always have $\text{Min A}(M,g) \leq \mathcal{A}_{\pi_2}(M,g)$. Gromov and Sormani in \cite{Sormani2023Conjectures} proposed the following conjecture.

\begin{conjecture}\label{conjecture Sormani and Gromov}
    Let $\alpha, D,V$ be positive numbers and consider $(M_j,g_j)_j$ a sequence of closed and  oriented $3$-dimensional Riemannian manifolds without boundary satisfying \[R_j \geq 0, \qquad \text{Min A}(M_j,g_j) \geq \alpha >0, \qquad \text{diam}(M_j) \leq D \qquad \text{and} \qquad \vol(M_j) \leq V.\] Then there exists a subsequence which converges in the volume preserving intrinsic flat sense to a three dimensional rectifiable geodesic metric space with Euclidean tangent cones satisfying generalized notions of nonnegative scalar curvature.
\end{conjecture}

The sequence $(M_j,g_j)_j$ constructed in Theorem \ref{tunnel counterexample} satisfies the hypothesis of scalar curvature, diameter and volume, but $\text{Min A}(M_j,g_j) \longrightarrow 0$ since $\text{Min A}(M_j,g_j)$ can be bounded above by the area of the $2$-spheres given by the slices of the tunnel, which goes to zero because the tunnel gets thinner. So this construction does not give a counterexample to the above conjecture. To the best of our knowledge, Conjecture \ref{conjecture Sormani and Gromov} remains open.

\begin{rmk}
    The same phenomenon occurs when we lift the construction of Theorem \ref{tunnel counterexample} to the universal cover. In this case, since we are dealing with noncompact manifolds, the right notion of convergence should be the pointed volume intrinsic flat. The universal covers provide a sequence of complete and oriented $3$-manifolds $(\widetilde{M}_j, \widetilde{g}_j)$ which are diffeomorphic to $S^2 \times \R$ and satisfy \[\mathcal A_{\pi_2}( \widetilde{M}_j, \widetilde{g}_j) \inf_{\widetilde{M}_j} R_{\widetilde{g}_j} \xrightarrow{\ j \rightarrow +\infty \ } 8\pi.\]  The upstairs picture should be $S^2\times \R$ with a bubble joined by a tunnel which appears near a fixed point in every fundamental region. The cylindrical part converges to $S^2\times\R$, while the $S^3$-pieces survive in this convergence. Therefore the pointed limit is \[M_\infty = (S^2\times\mathbb R) \sqcup \coprod_{k\in\mathbb Z}S^3\] and thus the pointed volume preserving intrinsic flat limit for this sequence is not the standard cylinder $\left(S^2\times\mathbb R, \ g_{S^2}+dt^2\right).$
\end{rmk}

\section{Rotationally Symmetric Stability}\label{Section warped}

Once the stability conjecture of Theorem \ref{BBN theorem} has been proven false, we show here and in the next section that adding certain hypotheses it is possible to make it true. First we restrict ourselves to Riemannian manifolds $S^2 \times S^1$ with rotationally symmetric metrics $g_\rho \coloneqq \rho(t)^2 g_{S^2}+dt^2$ for some smooth positive function $\rho$ defined on $S^1$. Since $S^3$ has $\pi_2(S^3) = 0$, we have to rule out the metrics that parameterize $3$-spheres. That is, since the model manifolds $(M,g)$ are closed and oriented $3$-dimensional Riemannian manifolds with $\pi_2(M) \neq 0$, we consider $(M,g)=(S^2 \times S^1,\ \rho(t)^2 g_{S^2}+dt^2)$ where $t\in\mathbb R/\ell\mathbb Z$ for some $\ell>0$ and $\rho$ is a positive smooth function. Geometrically $\ell$ is the length of the $S^1$-factor and when it is necessary to emphasize that, we will denote $S^1_\ell$. If it is convenient, we can interpret $\rho \colon \R \rightarrow (0,+\infty)$ as a function satisfying $\rho(x+\ell) = \rho(x)$ for all $x \in \R$. 

In the general case it was proved in \cite{MeeksYau1980, HassScott1988} that there exists a smooth immersion $f \colon S^2 \rightarrow M$ which is nontrivial in $\pi_2(M)$ and attains the infimum of $\mathcal{A}_{\pi_2}(M,g)$, that is, \[\text{Area}_{f^*g}(S^2) = \mathcal{A}_{\pi_2}(M,g).\] But in the case of a rotationally symmetric metric, Lemma \ref{minimizer slice} below ensures that it is possible to say more about the minimizers of $\mathcal{A}_{\pi_2}(M,g)$.

\begin{lemma}\label{minimizer slice}
    Consider $M = S^2 \times S^1$ with a rotationally symmetric metric $g_\rho \coloneqq  \rho^2(\theta)g_{S^2} + d\theta^2$. If $\theta_0$ is such that $\rho(\theta_0) = \min\limits_{\theta \in S^1} \rho(\theta) >0$ then the $2$-sphere $\Sigma_{\theta_0}\coloneq  S^2 \times \{\theta_0\}$ realizes $\mathcal{A}_{\pi_2}(M,g_\rho)$, that is, $\mathcal A_{\pi_2}(M,g_\rho) = \text{Area}_{g_{\rho}}(\Sigma_{\theta_0})$.
\end{lemma}

\begin{proof}
    Let $f \colon S^2 \rightarrow M$ be a smooth map which is nontrivial in $\pi_2(M)$ and write $f(v) = (u(v),\theta(v))$. Take local coordinates $(x_1,x_2)$ of $S^2$ and set for $i\in \{1,2\},$ \[u_i\coloneqq \partial_iu \qquad \text{and} \qquad \theta_i \coloneqq \partial_i\theta.\] Hence $\partial_if = (u_i,\theta_i)$ and the coefficients of the pullback tensor $f^*g$ are expressed as: \[(f^*g)_{ij} = \theta_i\theta_j + \rho(\theta)^2 \langle u_i,u_j\rangle_{S^2}.\] The determinant of the corresponding matrix is
    \begin{align*}
    \det\big((f^*g)_{ij}\big) &= \left(\theta_1^2 + \rho(\theta)^2\vert u_1\vert ^2 \right)\left(\theta_2^2 + \rho(\theta)^2\vert u_2\vert ^2 \right) - \left(\theta_1\theta_2 + \rho(\theta)^2\langle u_1,u_2\rangle_{S^2} \right)^2 \\
    &= \rho^2 (\theta)\underbrace{  \left(\theta_2^2\vert u_1\vert ^2 +\theta_1^2\vert u_2\vert ^2  -2 \theta_1\theta_2 \langle u_1,u_2\rangle_{S^2}\right)}_{(I)} + \rho^4(\theta) \underbrace{\left(\vert u_1\vert ^2\vert u_2\vert ^2 - \langle u_1,u_2\rangle_{S^2}^2 \right)}_{(II)}.
    \end{align*}
    Observe that $(I)$ and $(II)$ are nonnegative. In fact, $(II) \geq 0$ by the Cauchy-Schwarz inequality and $(I)  = \vert \theta_1u_2 - \theta_2 u_1\vert ^2 \geq 0.$ Moreover $(II) = \det\left((u^*g_{S^2})_{ij}\right)$ and thus we have \[\det((f^*g)_{ij}) \geq \rho(\theta)^4\det\left((u^*g_{S^2})_{ij}\right) \geq \rho(\theta_0)^4\det\left((u^*g_{S^2})_{ij}\right).\] 
    
    By hypothesis, $f$ corresponds to a nonzero class of $\pi_2(M)$. Hence $u \colon S^2 \rightarrow S^2$ has nonzero degree. This implies that  $\int_{S^2}\sqrt{\det\left((u^*g_{S^2})_{ij}\right)} \geq 4\pi \vert\text{deg}(u)\vert \geq 4\pi$. Thus,
    \begin{align*}
       \text{Area}_{f^*g}(S^2) &= \int_{S^2}\sqrt{\det((f^*g)_{ij})} \geq \int_{S^2} \rho(\theta_0)^2\sqrt{\det\left((u^*g_{S^2})_{ij}\right)} \geq 4\pi\rho(\theta_0)^2.
    \end{align*} Hence \[\text{Area}_{f^*g}(S^2) \geq 4\pi\rho(\theta_0)^2 = \text{Area}_{g_\rho}(S^2\times \{\theta_0\}).\] Taking the infimum over $f$, the result follows.
\end{proof}

Consider $M = S^2 \times S^1_\ell$ with rotationally symmetric metric $g_\rho$. When we analyze Theorem \ref{BBN theorem}, doing a rescale it is always possible to assume that $\inf_M R_{g_\rho} = 2.$ By Lemma \ref{minimizer slice}, take a minimizer slice $S^2 \times \{t_0\}$ where $\rho(t_0) = \min \rho$. Since $\rho$ is smooth and periodic, $\rho'(t_0) =0$ and $\rho''(t_0) \geq 0.$ Furthermore, after translating the coordinates if necessary, we can suppose that $t_0 = 0$ and if we denote $m := \rho(0)$, then $m = \min \rho$. We consider the almost equality condition \[\mathcal{A}_{\pi_2}(M,g_{\rho})= 4\pi m^2 \geq  4\pi(1-\varepsilon)\] for small $\varepsilon>0$. In this setting, $m$ satisfies the inequality $m^2\ge 1-\varepsilon$. Note also that if $\mathcal{A}_{\pi_2}(M,g_{\rho})\inf_M R_{g_\rho} \leq 8\pi$, then $m \leq 1$. And since $\rho$ is smooth, then for $t \in [-\ell,\ell]$, 
\begin{align} \label{DVM}
    \vert \rho(t)-m\vert = \vert \rho(t) - \rho(0)\vert  \leq \sup_{s \in [-\ell,\ell]}\vert \rho'(s)\vert  \vert t\vert  \leq \ell \ \sup\vert \rho'(s)\vert . 
\end{align}
Observe that, since the scalar curvature $R_{g_\rho}$ of a $3$-dimensional rotationally symmetric metric is given by \[R_{g_\rho} = -4 \dfrac{\rho''}{\rho} + 2\dfrac{1-(\rho')^2}{\rho^2},\] the assumption that $\inf_M R_{g_\rho} = 2$ gives 
\begin{equation}\label{Scal condition 1}
    \rho''\leq \dfrac{1-(\rho')^2-\rho^2}{2\rho} \leq \dfrac{1-\rho^2}{2\rho}.    
\end{equation}
But $m=\inf \rho >0$ and so 
\begin{equation}\label{rho and epsilon}
    1-\rho^2 \leq 1-m^2 \leq \varepsilon.    
\end{equation}
This implies that
\begin{equation}\label{Scal condition 2}
    \rho'' \leq \dfrac{\varepsilon}{2\rho}\leq \dfrac{\varepsilon}{2m} \leq \dfrac{\varepsilon}{2\sqrt{1-\varepsilon}}.  
\end{equation}
By the Fundamental Theorem of Calculus, we have for $0\leq t \leq \ell$, \[\rho'(t) = \int_{0}^t\rho''(s)ds \leq t \dfrac{\varepsilon}{2\sqrt{1-\varepsilon}} \leq  \dfrac{\ell\varepsilon}{2\sqrt{1-\varepsilon}}.\] And for $-\ell \leq t \leq 0$, \[-\rho'(t) = \int_t^0\rho''(s)ds \leq -t\dfrac{\varepsilon}{2\sqrt{1-\varepsilon}} \leq \dfrac{\ell\varepsilon}{2\sqrt{1-\varepsilon}}.\] Thus, if $0\leq t \leq \ell$ and $\rho'(t)\geq0$, then \[\vert \rho'(t)\vert  \leq \dfrac{\ell\varepsilon}{2\sqrt{1-\varepsilon}}.\] But if for this range of $t$ we have $\rho'(t) \leq 0$, then we can consider $\tilde{t} = t-\ell$. So $-\ell \leq \tilde{t} \leq 0$ and by the calculation above, \[- \rho'(\tilde{t} ) \leq \dfrac{\ell\varepsilon}{2\sqrt{1-\varepsilon}}.\] Since $\rho(t) = \rho(t+\ell)$ then $\rho'(\tilde{t}) = \rho'(t)$ and so \[\vert \rho'(t)\vert = -\rho'(t) \leq \dfrac{\ell\varepsilon}{2\sqrt{1-\varepsilon}}\] in this case. In the same fashion we conclude that for $-\ell \leq t \leq 0$ we also have $\vert \rho'(t)\vert  \leq \frac{\ell\varepsilon}{2\sqrt{1-\varepsilon}}.$ Hence by \eqref{DVM} it follows that for all $t\in [-\ell,\ell]$, \[\vert \rho(t) - m\vert  \leq \dfrac{\ell^2\varepsilon}{2\sqrt{1-\varepsilon}}.\] Therefore, as $\sqrt{1-\varepsilon} \leq m \leq 1$, \[\rho(t) - 1 \leq \dfrac{\ell^2\varepsilon}{2\sqrt{1-\varepsilon}}  + m -1 \leq \dfrac{\ell^2\varepsilon}{2\sqrt{1-\varepsilon}} \] and \[1-\rho(t) \leq 1 -m+ \dfrac{\ell^2\varepsilon}{2\sqrt{1-\varepsilon}} \leq \varepsilon + \dfrac{\ell^2\varepsilon}{2\sqrt{1-\varepsilon}}\] where the last inequality holds if $0<\varepsilon < 1$. Hence
\begin{equation*}
\vert \rho(t)-1\vert  \leq \varepsilon + \dfrac{\ell^2\varepsilon}{2\sqrt{1-\varepsilon}}    
\end{equation*}
for all $t \in [-\ell,\ell]$. Notice that we always have the bound $\ell \leq 2\text{diam}_{g_\rho}(M)$ and hence we can write the inequalities above as
\begin{equation}\label{estimativa de f}
\vert \rho(t)-1\vert  \leq \varepsilon + \dfrac{2\varepsilon}{\sqrt{1-\varepsilon}}\ \text{diam}_{g_\rho}(M)^2 \qquad \text{and} \qquad \vert \rho'(t)\vert  \leq \frac{\varepsilon}{\sqrt{1-\varepsilon}}\ \text{diam}_{g_\rho}(M)
\end{equation}
which, by periodicity, hold for every $t \in \R$. 

We have then proved the following bound for the $C^1$-norm of the warping function. 

\begin{teo}\label{teorema C1}
    Given $L>0$ and $\delta >0$, there exists $\varepsilon = \varepsilon(L,\delta) \in (0,1)$ such that for every closed and oriented manifold $M = S^2 \times S^1$ with rotationally symmetric metric $g_{\rho}$ which satisfies
    \begin{enumerate}
        \item[1)] diam($M$) $\leq L$;
        \item[2)] $\inf_M R_{g_\rho} =2$;
        \item[3)] $\mathcal{A}_{\pi_2}(M,g_{\rho}) \geq 4\pi(1-\varepsilon)$,
    \end{enumerate}
    its warping function is $\delta$ close to the constant function $1$ in the $C^1$-norm, that is, $ \Vert \rho-1\Vert_{C^1} \leq \delta$.
\end{teo}

Theorem \ref{teorema C1} proves a $C^1$-type of stability result for Theorem \ref{BBN theorem} in the context of rotationally symmetric manifolds.

\begin{rmk}
    The fact that the warping function is $C^1$-close to the constant function equals to one is not enough to ensure a $C^1$-convergence when we are in the context of a sequence of rotationally symmetric manifolds. In fact, if $M_j = S^2 \times S^1_{\ell_j}$ is a sequence satisfying the hypotheses of Theorem \ref{teorema C1} with $\varepsilon_j \rightarrow0$, it may happen that the $S^1$-factor $\ell_j$ of each $M_j$ collapses. But if $\ell_j \rightarrow \ell >0$, after a identification of the circles $S^1_{\ell_j}$  with $S^1_\ell$ we have the $C^1$-convergence of $(S^2 \times S^1_{\ell_j}, g_{\rho_j})$ to $(S^2\times S^1_\ell, g_{S^2}+dt^2)$ as expected.
\end{rmk}

Theorem \ref{teorema C1} allows to prove Theorem \ref{Theorem C} in two parts: first we check in Subsection \ref{Subsection GH} that it implies Gromov-Hausdorff convergence in the context of sequences. After that, in Subsection \ref{Subsection IF} we guarantee that the Gromov-Hausdorff convergent subsequence, up to a further subsequence, also converges in the volume preserving intrinsic flat sense. Thus the proof of Theorem \ref{Theorem C} follows from a combination of Corollary \ref{GH convergence warped} and Corollary \ref{SWIF convergence warped} below.

\subsection{Gromov-Hausdorff Convergence in the Rotationally Symmetric Case}\label{Subsection GH}

The proof of Theorem \ref{teorema C1} also gives an explicit estimate in terms of the Gromov-Hausdorff distance. For a manifold $M=S^2\times S^1_\ell$ with rotationally symmetric metric $g_\rho$, by $(\ref{estimativa de f})$ we have the inequality:
\begin{equation}\label{C0 estimate of warping function}
    \Vert \rho-1\Vert_{C^0} \leq \varepsilon + \dfrac{2\varepsilon}{\sqrt{1-\varepsilon}} \text{diam}_{g_\rho}(M)^2 =: \eta(\rho).
\end{equation}
Fixing the standard product metric $g_1 \coloneq g_{S^2}+dt^2$ in  $M = S^2\times S^1_\ell$, this implies that if $\varepsilon$ is small enough the distances $d_{g_\rho}$ and $d_{g_1}$ are related by \[(1-\eta)d_{g_1} \leq  d_{g_\rho} \leq (1+\eta)d_{g_1}.\] Hence \[\vert d_{g_\rho}-d_{g_1}\vert \leq \eta \  d_{g_1}.\] By Proposition \ref{GH isometry} it is possible to conclude that
\begin{equation}\label{GH estimate}
    d_{GH}((M,g_\rho),(M,g_1)) \leq 2\eta \ \text{diam}_{g_1}(M) = 2\left(\varepsilon + 2\dfrac{\text{diam}_{g_\rho}(M)^2\varepsilon}{\sqrt{1-\varepsilon}}\right) \text{diam}_{g_1}(M).
\end{equation}

In terms of sequences, we expect that such an estimate will allow convergence, at least for a subsequence. For $M_j = (S^2\times S^1_{\ell_j},g_{\rho_j})$ with $\mathcal{A}_{\pi_2}(M_j,g_{\rho_j}) \geq 4\pi(1-\varepsilon_j)$ and $\varepsilon_j \rightarrow 0$, the length $\ell_j \in (0,2L]$ if $\text{diam}(M_j) \leq L$. Then there exists a subsequence that converges to some number $\ell_\infty \in [0,2L]$. In this case the Gromov-Hausdorff limit should be $(S^2 \times  S^1_{\ell_\infty},g_1)$. But if $\ell_\infty = 0$, then the $S^1$-factor collapses and the limit is not $3$-dimensional anymore. To avoid such thing, we can impose the condition $\text{Vol}(M_j) \geq v >0$ for all index $j$.

\begin{cor}\label{GH convergence warped}
    Let $M_j = S^2 \times S^1_{\ell_j}$ be a sequence of closed oriented manifolds with rotationally symmetric metric $g_{\rho_j}$ such that:
    \begin{enumerate}
        \item[1)] $\text{diam}_{g_{\rho_j}}(M_j) \leq \ L$ \ \  and \ \  $\vol_{g_{\rho_j}}(M_j) \geq v >0$;
        \item[2)] $\inf_{M_j} R_{g_{\rho_j}} =2$;
        \item[3)] $\mathcal{A}_{\pi_2}(M_j,g_{\rho_j}) \geq 4\pi(1-\varepsilon_j)$ with $\varepsilon_j \rightarrow 0$,
    \end{enumerate}
    where $L$ and $v$ are positive constants that do not depend on $j$ and $\varepsilon_j \in (0,1)$. Then there exists a subsequence $(M_{j_k})_k$ that converges in the Gromov-Hausdorff sense to $S^2\times S^1_{\ell_\infty}$, $\ell_\infty \neq 0$, with canonical product metric $g_1 = g_{S^2}+dt^2$.
\end{cor}
\begin{proof}
    Fix an index $j\in \mathbb N$. In the proof of Theorem \ref{teorema C1}, estimate \eqref{rho and epsilon} shows that $\sqrt{1-\varepsilon_j} \leq \rho_j$. By inequality (\ref{C0 estimate of warping function}) we have an upper bound $\rho_j\leq 1+\eta_j$ where $\eta_j \coloneqq \varepsilon_j+\frac{2\varepsilon_j}{\sqrt{1-\varepsilon_j}}\ \text{diam}_{g_{\rho_j}}(M_j)^2$. Hence, when comparing to the Riemannian manifold $(S^2\times S^1_{\ell_j},g_1)$, we have \[\sqrt{1-\varepsilon_j} \ d_{g_1} \leq d_{g_{\rho_j}} \leq (1+\eta_j)d_{g_1}\] and so \[d_{g_{\rho_j}}-d_{g_1} \leq \eta_j\ d_{g_1} \qquad \text{and} \qquad d_{g_1}-d_{g_{\rho_j}} \leq (1-\sqrt{1-\varepsilon_j})\ d_{g_1}.\] As $1-\sqrt{1-\varepsilon_j} = \frac{\varepsilon_j}{1+\sqrt{1-\varepsilon_j}} \leq \varepsilon_j \leq \eta_j$, then \[\vert d_{g_{\rho_j}}-d_{g_1}\vert \leq \eta_j \ d_{g_1}.\] And note that as $\varepsilon_j \rightarrow 0$, we also have $\eta_j \rightarrow 0$ since the diameter of the manifolds is uniformly bounded. 
    
    For a subsequence $j_k$ we have $\ell_{j_k}\xrightarrow{k\rightarrow +\infty} \ell_\infty$ for some $\ell_\infty \in[0,2L].$ We claim that $\ell_\infty$ is not zero. In fact, \[0<v \leq \vol_{g_{\rho_{j_k}}}(M_{j_k})= 4\pi \int_0^{\ell_{j_k}}\rho_{j_k}^2(t) \ dt \leq 4\pi(1+\eta_{j_k})^2\ell_{j_k}\] and since $\eta_{j_k}$ is uniformly bounded by hypothesis $1)$, it follows that $\ell_\infty >0.$ By (\ref{GH estimate}), if $M = S^2 \times S^1_{\ell_\infty}$,
    \begin{align*}
        d_{GH}\big((M_{j_k},g_{\rho_{j_k}}), (M,g_1)\big) &\leq
    \begin{aligned}[t]
    d_{GH}\big((S^2\times S^1_{\ell_{j_k}},g_{\rho_{j_k}})&,
    (S^2\times S^1_{\ell_{j_k}},g_1)\big)\\
    &+d_{GH}\big((S^2\times S^1_{\ell_{j_k}},g_1), (S^2\times S^1_{\ell_\infty},g_1)\big)
    \end{aligned}\\
    &\leq 2\eta_{j_k} \ \text{diam}_{g_1}(S^2\times S^1_{\ell_{j_k}}) +d_{GH}\big((S^2\times S^1_{\ell_{j_k}},g_1), (S^2\times S^1_{\ell_\infty},g_1)\big).
    \end{align*} 
and the right hand side converges to zero as $k \rightarrow+\infty$.
\end{proof}

\subsection{Intrinsic Flat Convergence in the Rotationally Symmetric Case}\label{Subsection IF}

Consider again a sequence of manifolds $M_j=S^2\times S^1_{\ell_j}$ with rotationally symmetric metric $g_{\rho_j}$ satisfying $\inf_{M_j} R_{g_{\rho_j}} =2$ and $\mathcal A_{\pi_2}(M_j,g_{\rho_j}) \geq 4\pi(1-\varepsilon_j)$ with $\varepsilon_j \rightarrow0$. In a warped product metric the Ricci curvature is given as follows: if $\partial_t$ denotes the unit radial direction and $X$ is any unit vector tangent to the $S^2$-factor, \[\Ric_{g_{\rho_j}}(X,X) = \dfrac{1-(\rho_j')^2}{\rho_j^2} - \dfrac{\rho_j''}{\rho_j}, \qquad \Ric_{g_{\rho_j}}(\partial_t, \partial_t)= -2\dfrac{\rho_j''}{\rho_j} \qquad \text{and} \qquad \Ric_{g_{\rho_j}}(\partial_t,X) =0.\] If all the manifolds have uniformly bounded diameter, that is, if for some $L>0$, $\text{diam}_{g_{\rho_j}}(M_j) \leq L$ for all $j\in \mathbb{N}$, the inequality \eqref{estimativa de f} gives: \[\Vert \rho_j-1 \Vert_{C^0} \leq \varepsilon_j + \dfrac{2\varepsilon_j L^2}{\sqrt{1-\varepsilon_j}} =: \delta_0^j \qquad \text{and} \qquad \Vert \rho_j' \Vert_{C^0} \leq \dfrac{\varepsilon_j L}{\sqrt{1-\varepsilon_j}} =:\delta_1^j.\] Also the scalar condition $R_{g_{\rho_j}} \geq 2$ gives, by equations \eqref{Scal condition 1}  and \eqref{Scal condition 2}, that for sufficiently large $j$ the second derivative of the warping function satisfies $\rho_j'' \leq \varepsilon_j.$ 

After discarding finitely many terms, we can always suppose that $\varepsilon_j\leq \frac{1}{2}$ for all $j \in \mathbb{N}$. Hence $\frac{1}{2} \leq 1-\varepsilon_j \leq 1$. In the proof of Theorem \ref{teorema C1}, estimate \eqref{rho and epsilon} shows that $\sqrt{1-\varepsilon_j} \leq \rho_j$. This implies that $\frac{1}{\sqrt{2}}\leq \sqrt{1-\varepsilon_j}\leq\rho_j$ and so $\frac{1}{\rho_j}\leq \sqrt{2}$ and $\frac{1}{\rho_j^2} \leq 2.$ With this we have \[\dfrac{1-(\rho_j')^2}{\rho_j^2} \geq 0 -\dfrac{(\delta_1^j)^2}{\rho_j^2} \geq  -2(\delta_1^j)^2 \qquad \text{and} \qquad -\dfrac{\rho_j''}{\rho_j} \geq -\dfrac{\varepsilon_j}{\rho_j} \geq -\sqrt{2}\varepsilon_j.\] Combining this with the expression of Ricci curvature and the fact that $\varepsilon_j \in (0,\frac{1}{2})$ we conclude that $\Ric_{g_{\rho_j}} \geq \Lambda$ where $\Lambda$ is a constant that does not depend on $j$. 

Also note that since $0 \leq \rho_j \leq 1+\delta_0^j$, then \[\vol_{g_{\rho_j}}(M_j)= 4\pi \int_0^{\ell_j}\rho_j^2(t) \ dt \leq (1+\delta_0^j)^2 4\pi \int_0^{\ell_j}dt =  (1+\delta_0^j)^2\vol_{g_1}(M_j)\] and thus \[\vol_{g_{\rho_j}}(M_j) - \vol_{g_1}(M_j) \leq ((1+\delta_0^j)^2-1)\vol_{g_1}(M_j).\] On the other hand, since $1-\varepsilon_j \leq \rho_j^2$, then \[(1-\varepsilon_j)\vol_{g_1}(M_j) \leq \vol_{g_{\rho_j}}(M_j)\] and so \[\vol_{g_1}(M_j) -\vol_{g_{\rho_j}}(M_j) \leq \varepsilon_j \vol_{g_1}(M_j).\] Hence it holds that  \[\vert\vol_{g_{\rho_j}}(M_j)-\vol_{g_1}(M_j)\vert \leq \max\left\{(1+\delta_0^j)^2-1, \varepsilon_j\right\}\vol_{g_1}(M_j).\] Notice that $\ell_j\leq 2\text{diam}_{g_{\rho_j}}(M_j) \leq 2L$ implies $\vol_{g_1}(M_j) \leq 8\pi L$. Therefore, by the fact that $\delta_0^j \rightarrow 0$ as $\varepsilon_j \rightarrow 0$ it follows that 
\begin{equation}\label{convergence of volumes}
\vert\vol_{g_{\rho_j}}(M_j)-\vol_{g_1}(M_j)\vert \xrightarrow{\ j\rightarrow +\infty \ } 0.    
\end{equation}

Now it can be proved the volume preserving intrinsic flat convergence in the rotationally symmetric case.

\begin{cor} \label{SWIF convergence warped}
    Let $M_j = S^2 \times S^1_{\ell_j}$ be a sequence of closed oriented manifolds with rotationally symmetric metric $g_{\rho_j}$ such that:
    \begin{enumerate}
        \item[1)] $\text{diam}_{g_{\rho_j}}(M_j) \leq \ L$ \ \ \text{and} \ \  $\vol_{g_{\rho_j}}(M_j) \geq v >0$;
        \item[2)] $ \inf_{M_j} R_{g_{\rho_j}} =2$;
        \item[3)] $\mathcal{A}_{\pi_2}(M_j,g_{\rho_j}) \geq 4\pi(1-\varepsilon_j)$ with $\varepsilon_j \rightarrow 0$,
    \end{enumerate}
    where $L$ and $v$ are positive constants that do not depend on $j$ and $\varepsilon_j \in (0,1)$. Then there exists a subsequence $(M_{j_k})_k$ that converges to $M_1 =S^2\times S^1_{\ell_\infty}$, $\ell_\infty \neq 0$, with canonical product metric $g_1 = g_{S^2}+dt^2$, in the volume preserving intrinsic flat sense, that is, \[M_{j_k} \xrightarrow{\ \mathcal{F} \ }{M_1}\qquad \text{and} \qquad \vol(M_{j_k}) \xrightarrow{\ k \rightarrow +\infty\ } \vol(M_1).\]
\end{cor}

\begin{proof}
    By Corollary \ref{GH convergence warped} there exists a subsequence $(M_{j_k})_k$ that converges in the Gromov-Hausdorff sense to $M_1 = S^2 \times S^1_{\ell_\infty}$ with the standard metric. The triangle inequality implies that the volumes satisfy the relation \[\vert \vol_{g_{\rho_{j_k}}}(M_{j_k}) - \vol_{g_1}(M_1)\vert \leq \vert \vol_{g_{\rho_{j_k}}}(M_{j_k}) - \vol_{g_1}(M_{j_k})\vert + \vert \vol_{g_1}(M_{j_k}) - \vol_{g_1}(M_1)\vert \] and the convergence of the volumes \eqref{convergence of volumes} ensures that the first summand of the right hand side goes to zero. The second term also goes to zero since \[\vol_{g_1}(M_{j_k}) = 4\pi \ell_{j_k} \longrightarrow 4\pi \ell_\infty = \vol_{g_1}(M_1).\] Therefore \[\vert \vol_{g_{\rho_{j_k}}}(M_{j_k}) - \vol_{g_1}(M_1)\vert \longrightarrow 0.\] Identifying these Riemannian manifolds as currents, the convergence of the volumes translates to the convergence of the masses. Along with hypothesis $1)$, Wenger's compactness theorem \ref{Wenger compactness theorem} guarantees a further subsequence converging in the intrinsic flat sense to a subset of $M_1$. Since the Ricci curvature is uniformly bounded below, Theorem A in \cite{MatveevPortegies2017} identifies this limit with $M_1$. 
\end{proof}

\section{$\Ric \geq 0$ Stability}\label{Section Ric}

The main purpose of this section is to prove the stability of Theorem \ref{BBN theorem} under the additional hypothesis of $\Ric \geq 0$. Since the compactness theorems of Gromov-Hausdorff distance work under the hypothesis of a lower bound on the curvature, to expect such a convergence when there is no restriction in the geometry of the Riemannian manifolds it is reasonable to add the condition of a lower bound for Ricci. In view of the fact that $S^2 \times \R$ is not compact, it is expected to obtain a pointed Gromov-Hausdorff convergence of the universal covers. An important ingredient in the proof of Theorem \ref{Ricci theorem} is the topological gap theorem of Xu \cite{Xu2023TopologicalGap}, which is written below.

\begin{teo}[Theorem 1.1 in \cite{Xu2023TopologicalGap}]\label{Theorem of Xu}
    Suppose $M$ is a closed and oriented $3$-manifold with $\pi_2(M) \neq 0$. If $M$ is not covered by $S^2\times S^1$, then for any metric $g$ on $M$ we have \[\mathcal A_{\pi_2}(M,g) \inf_M R_g \leq 24\pi \cdot \dfrac{2-\sqrt{2}}{4-\sqrt{2}} \ \ (\approx 5.44\pi).\]
\end{teo}

Therefore, for a sequence $(M_j,g_j)_j$ satisfying the almost rigidity hypothesis \[\mathcal A_{\pi_2}(M_j,g_j)\inf_{M_j}R_{g_j}\xrightarrow{\ j \rightarrow +\infty \ } 8\pi,\] Theorem \ref{Theorem of Xu} implies that, for sufficiently large $j$, $M_j$ must be a quotient of $S^2\times S^1$. This provides a topological reduction in the argument and allows us to restrict the subsequent analysis to manifolds with the topology of the rigidity model in the Bray-Brendle-Neves theorem.

\begin{proof}[Proof of Theorem \ref{Ricci theorem}]
    As $\varepsilon_j \rightarrow0$, we can suppose that  $8\pi(1-\varepsilon_j) > 24\pi \cdot \frac{2-\sqrt{2}}{4-\sqrt{2}}$. By Theorem \ref{Theorem of Xu} each $M_j$ is covered by $S^2 \times S^1$. This implies that the universal cover $\widetilde{M_j}$ of $M_j$ is $S^2\times \R$, which has two ends and, therefore, a line. As $\Ric_{g_j} \geq 0$, Cheeger-Gromoll splitting theorem guarantees that $\widetilde{M_j}$, with the cover metric $\widetilde{g}_j$, is isometric to a Riemannian product $(N_j \times\R, h_j+dt^2)$ where $N_j$ is a $2$-sphere. Furthermore, the sectional curvature  $K_{h_j}$ of $(N_j,h_j)$ satisfies \[K_{h_j} = \dfrac{1}{2}R_{h_j} = \dfrac{1}{2}R_{g_j} \geq 1. \] By the Gauss-Bonnet theorem, \[4\pi = 2\pi \mathcal{X}(S^2) = \int_{N_j} K_{h_j} \geq \text{Area}_{h_j}(N_j).\] However, \[\text{Area}_{h_j}(N_j) \geq \mathcal{A}_{\pi_2}(\widetilde{M}_j, \tilde{g}_j) = \mathcal{A}_{\pi_2}(M_j,g_j) \geq 4\pi(1-\varepsilon_j)\] and so \[4\pi(1-\varepsilon_j) \leq \text{Area}_{h_j}(N_j) \leq 4\pi\] which implies, as $\varepsilon_j \rightarrow0$, that $\text{Area}_{h_j}(N_j) \rightarrow4\pi.$ 
    
    Also, the condition $K_{h_j} \geq 1$ implies, by Bonnet-Myers theorem, that $\text{diam}(N_j) \leq \pi$. Proposition \ref{Gromov compactness} ensures that the class of $2$-dimensional Riemannian manifolds with uniformly bounded diameter and  uniform lower sectional curvature bound is pre-compact in the Gromov-Hausdorff topology. Hence, after passing to a subsequence and relabeling, we may assume that there exists an Alexandrov space $N$ with Hausdorff dimension $\text{dim}_H(N) \leq 2$ such that \[N_j \xrightarrow{\ GH \ } N\] and its sectional curvature (in the Alexandrov sense) is $K \geq 1$. Moreover, by Theorem \ref{Bonnet Myers Alexandrov}, $\text{diam}(N) \leq \pi.$
    
    For all $j \in \mathbb{N}$ we can use Theorem \ref{Bishop gromov equality alexandrov} to conclude that, for all $p_j \in N_j$ and $r>0$, \[\mathcal H^2(B_r(p_j)) \leq V_r^1\] where $\mathcal H^2(B_r(p_j))$ is the $2$-dimensional Hausdorff measure of the ball of radius $r>0$ centered at $p_j$ and $V_r^1$ is the volume of the ball of radius $r$ in the round sphere $S^2$ of curvature $1$. Regarding the lower bound, since we have diam$(N_j) \leq \pi$, Theorem \ref{Bishop Gromov alexandrov} implies that, for all $0<r \leq \pi$, \[ \dfrac{\mathcal H^2(B_r(p_j))}{V_r^1} \geq \dfrac{\mathcal H^2(B_\pi(p_j))}{V_\pi^1} = \dfrac{\text{Area}_{h_j}(N_j)}{4\pi} \geq \dfrac{4\pi(1-\varepsilon_j)}{4\pi} = 1-\varepsilon_j.\] It also holds for $R>\pi$ that \[ \dfrac{\mathcal H^2(B_R(p_j))}{V_R^1} = \dfrac{\text{Area}_{h_j}(N_j)}{4\pi} \geq \dfrac{4\pi(1-\varepsilon_j)}{4\pi} = 1-\varepsilon_j.\] Therefore, $\mathcal H^2(B_r(p_j)) \geq(1-\varepsilon_j)V_r^1$ for all $r>0.$ Thus, given $r>0$, for all $j \in \mathbb N$,
    \begin{equation}\label{volume estimate BG}
        (1-\varepsilon_j)V_r^1 \leq \mathcal H^2(B_r(p_j)) \leq V_r^1.   
    \end{equation}
    Given $p \in N$, by Proposition \ref{pointed GH and GH COMPACT} take $p_j \in N_j$ such that $(N_j,p_j)\xrightarrow{\ pGH\ }{(N,p)}$. We would like to pass to the limit in the inequalities \eqref{volume estimate BG}. 
    
    In order to do that, first we show that the Hausdorff dimension of $N$ is exactly $2$, that is, the Gromov-Hausdorff limit of $(N_j)_j$ does not collapse. In fact, by the pointed Gromov-Hausdorff convergence and the fact that the diameters are uniformly bounded, there exists a sequence $\nu_j \rightarrow 0$ and $\nu_j$-isometries $f_j \colon N \rightarrow N_j$, that is, maps such that for all $x,y \in N$, \[\vert d_{h_j}(f_j(x),f_j(y)) - d_N(x,y)\vert \leq \nu_j \qquad \text{and} \qquad N_j \subseteq U_{\nu_j}(f_j(N)) \] where $U_{\nu_j}$ denotes the $\nu_j$-neighborhood. By adjusting the approximations and enlarging $\nu_j$ if necessary, we may also assume that $f_j(p)=p_j$. 
     
    If $\dim_H(N) <2$, then $\mathcal H^2(N) = 0.$ By the definition of Hausdorff measure and the compactness of $N$, for every $\delta>0$ there exists a finite cover by open balls \[N = \bigcup_{i=1}^m B_{r_i}(q_i)\qquad \text{with}  \qquad \sum\limits_{i=1}^m r_i^2 < \delta.\] We prove first that for each $j \in \mathbb N$, \[N_j \subseteq \bigcup_{i=1}^m B_{r_i + 2\nu_j}(f_j(q_i)).\] Indeed, given $y\in N_j$, there exists $x\in N$ such that $d_{h_j}(y,f_j(x)) < \nu_j$. Since the balls cover $N$, there exists $i \in \{1,\dots,m\}$ such that $x \in B_{r_i}(q_i)$. Therefore \[d_{h_j}(y,f_j(q_i)) \leq d_{h_j}(y,f_j(x)) + d_{h_j}(f_j(x),f_j(q_i)) < \nu_j + d_N(x,q_i) + \nu_j < r_i +2\nu_j.\] Theorem \ref{Bishop gromov equality alexandrov} and the inclusion above imply \[\mathcal H^2(N_j) \leq \sum\limits_{i=1}^m \mathcal H^2 \big(B_{r_i + 2\nu_j}(f_j(q_i))\big) \leq \sum\limits_{i=1}^m V_{r_i+2\nu_j}^1 \leq 2\pi \sum\limits_{i=1}^m (r_i+2\nu_j)^2.\] Letting $j\rightarrow +\infty$, \[\limsup\limits_{j\rightarrow +\infty} \mathcal H^2(N_j) \leq 2\pi \sum\limits_{i=1}^m r_i^2  < 2\pi\delta.\] Letting $\delta\rightarrow 0,$ it follows $\lim\limits_{j\rightarrow +\infty} \mathcal H^2(N_j) = 0.$ On the other hand, if $r_0 \in (0,\pi)$ by \eqref{volume estimate BG}, \[\mathcal H^2(N_j) \geq \mathcal H^2(B_{r_0}(p_j)) \geq (1-\varepsilon_j)V_{r_0}^1\] and so \[\liminf\limits_{j \rightarrow +\infty } \mathcal H^2(N_j) \geq V_{r_0}^1 >0\] which is a contradiction. Thus we conclude that $\dim_H(N)=2$.
    
    We now claim that $\mathcal H^2(B_r(p)) =V_r^1$ for every $r>0$, that is, \eqref{volume estimate BG} holds in the limit. By Theorem 10.8 in \cite{BuragoGromovPerelman1992}, for all sets $A \subseteq N$ with $\mathcal H^2(\partial A) = 0$, we have
    \begin{equation}\label{teorema 10.8}
        \mathcal H^2(U_{\nu_j} (f_j(A))) \xrightarrow{\ j \rightarrow +\infty \ }{} \mathcal H^2(A).    
    \end{equation}
    Notice that for almost all $r>0$, $\mathcal H^2(\partial B_r(p)) =0$ and so \eqref{teorema 10.8} holds for such $A= B_r(p)$.
    Consider $r> 2\nu_j$. The following chain of inclusions is true 
    \begin{equation}\label{chain of inclusions}
        B_{r-2\nu_j}(p_j) \subseteq U_{\nu_j}\big (f_j(B_r(p))\big) \subseteq B_{r+2\nu_j}(p_j).    
    \end{equation}
    In fact, for the first inclusion, take $y \in B_{r-2\nu_j}(p_j).$ Then there exists $x\in N$ such that $d_{h_j}(y,f_j(x)) < \nu_j$. Hence \[d_N(x,p) \leq d_{h_j}(f_j(x),f_j(p)) + \nu_j \leq d_{h_j}(f_j(x),y)+ d_{h_j}(y,p_j)+\nu_j < r.\] For the other inclusion, take $y \in U_{\nu_j} \big(f_j(B_r(p))\big)$. Then there exists $x \in B_r(p)$ such that $d_{h_j}(y,f_j(x)) < \nu_j$ and thus \[d_{h_j}(y,p_j) \leq d_{h_j}(y,f_j(x)) + d_{h_j}(f_j(x), f_j(p)) < \nu_j + d_N(x,p)+ \nu_j < r+2\nu_j.\]

    Now fix $r>0$ such that $\mathcal H^2(\partial B_r(p))=0$ and $0<\delta < r.$ Since $\nu_j \rightarrow 0$, for sufficiently large $j$ the chain of inclusions \eqref{chain of inclusions} implies \[\mathcal H^2(B_{r-\delta}(p_j)) \leq \mathcal H^2\big(U_{\nu_j}\big(f_j(B_r(p))\big)  \big) \leq \mathcal H^2(B_{r+\delta}(p_j)).\] By the volume estimate \eqref{volume estimate BG},
    \[(1-\varepsilon_j)V_{r-\delta}^1 \leq  \mathcal H^2\big(U_{\nu_j}\big(f_j(B_r(p))\big)  \big) \leq V_{r+\delta}^1.\] Letting $j\rightarrow +\infty$, by \eqref{teorema 10.8} it follows that \[V_{r-\delta}^1 \leq \mathcal H^2(B_r(p)) \leq V_{r+\delta}^1.\] Letting $\delta \rightarrow 0$, we conclude that $\mathcal H^2(B_r(p)) = V_r^1.$ For an arbitrary $r>0$, choose a sequence of radii $r_k\uparrow r$ with $\mathcal H^2(\partial B_{r_k}(p)) = 0.$ Hence \[\mathcal H^2(B_r(p)) = \lim\limits_{r_k \uparrow r} \mathcal H^2(B_{r_k}(p)) = \lim\limits_{r_k \uparrow r}  V_{r_k}^1 = V_r^1.\]
    
    Once $\mathcal H^2(B_r(p)) = V_r^1$ for all $r>0$, the equality case in Theorem \ref{Bishop gromov equality alexandrov} (see also Theorem 10.2 in \cite{BuragoGromovPerelman1992}) guarantees that $N$ is isometric to the round $2$-sphere of curvature $1$. This means that the Riemannian manifolds $(N_j,h_j)$ converge in the Gromov-Hausdorff sense to $(S^2,g_{S^2})$.
    
    Finally, given $x=(p,t)\in S^2\times \R$, take $p_j \in N_j$ as before, that is, such that $(N_j, p_j) \xrightarrow{\ pGH \ } (N,p)$. For the product metric $g \coloneqq  g_{S^2}+dt^2$ in $S^2\times \R$ we have that \[(\widetilde{M_j},(p_j,t)) \xrightarrow{\ pGH \ }{(S^2\times \R,(p,t)).}\] To prove that, we check properties (1), (2) and (3) of Definition \ref{definition pGH}. Indeed, the distances satisfy \[d_{\widetilde{g}_j}\left((x,s),(y,u)\right) = \sqrt{d_{h_j}(x,y)^2+\vert s-u\vert^2}\] for $(x,s),(y,u) \in N_j\times \R$ and \[d_{g}\left((x,s),(y,u)\right) = \sqrt{d_{S^2}(x,y)^2+\vert s-u\vert^2}\] for $(x,s),(y,u) \in S^2\times \R$. Given $r>0$ and $\varepsilon>0$, since $(N_j, p_j) \xrightarrow{\ pGH \ }(S^2,p)$ then there exists $j_0\in\mathbb{N}$ such that for all $j \geq j_0$ we have maps $f_j \colon B_r(p_j) \subset N_j \rightarrow S^2$ satisfying $(1)$,$(2)$ and $(3)$ in Definition \ref{definition pGH} for this $r>0$ and $\varepsilon' =\frac{\varepsilon}{4}$. Define
\[\begin{aligned}
    F_j \colon B_r((p_j,t)) \subset &\ N_j \times \R \rightarrow S^2 \times \R\\
    (x,s) & \longmapsto (f_j(x),s).
  \end{aligned}\]
It is easy to see that $F_j$ is well defined. In fact, if $(x,s)\in B_r((p_j,t))\subset N_j \times \R$ then $d_{h_j}(x,p_j) \leq d_{\widetilde{g}_j} ((x,s),(p_j,t)) <r$. Hence $x \in B_r(p_j)\subset N_j$.

Note that:

\begin{enumerate}
    \item[(1)] $F_j(p_j,t) = (f_j(p_j),t) = (p,t)$.
    \item[(2)] For $(x,s),(y,u) \in B_r((p_j,t)) \subset N_j \times \R$, if $a:= d_{h_j}(x,y)$, $b:= d_{S^2}(f_j(x),f_j(y))$ and $c:= \vert s-u\vert$, then: \[d_{\tilde{g}_j}\left((x,s),(y,u) \right) = \sqrt{a^2+c^2}\] and \[d_{g}\left(F_j(x,s),F_j(y,u) \right) = \sqrt{b^2+c^2}.\] Therefore \[\left\vert d_g\left(F_j(x,s),F_j(y,u) \right) -d_{\tilde{g}_j}\left((x,s),(y,u) \right) \right\vert \leq \vert a-b\vert \leq \text{dis }f_j <\varepsilon' < \varepsilon. \]
    \item[(3)] We show that the $\varepsilon$-neighborhood of $F_j\left(B_r\left((p_j,t)\right) \right)$ contains the ball $B_{r-\varepsilon}((p,t)) \subset S^2 \times \R.$ Take $(z,s) \in B_{r-\varepsilon}((p,t))\subset S^2 \times \R$. Put \[A := d_{S^2}(z,p) \qquad \text{and} \qquad B:= \vert s-t\vert.\] Then \[d_{g}\left((z,s),(p,t) \right)^2 = A^2+B^2 < (r-\varepsilon)^2.\] In particular, $A<r-\varepsilon < r-\varepsilon'$. Hence $z \in B_{r-\varepsilon'}(p)\subset S^2$. By property $(3)$ of $f_j$ in Definition \ref{definition pGH},  there exists $x \in B_r(p_j)\subset N_j$ such that $d_{S^2}(f_j(x),z) <\varepsilon'.$ As $\text{dis } f_j < \varepsilon'$ we have \[\vert d_{S^2}(f_j(x),p) - d_{h_j}(x,p_j)\vert = \vert d_{S^2}(f_j(x),f_j(p_j)) - d_{h_j}(x,p_j)\vert <\varepsilon'.\] Furthermore, \[d_{S^2}(f_j(x),p) \leq d_{S^2}(f_j(x),z) + d_{S^2}(z,p) < \varepsilon'+A.\] Thus \[d_{h_j}(x,p_j)< d_{S^2}(f_j(x),p) + \varepsilon' < 2\varepsilon'+A.\] With this it is easy to see that  $(x,s) \in B_r((p_j,t))\subset N_j \times \R$. In fact,
    \begin{align*}
        d_{\tilde{g}_j}\left((x,s),(p_j,t) \right)^2 &= d_{h_j}(x,p_j)^2+B^2 \leq (A+2\varepsilon')^2+B^2=(A^2+B^2) + 4\varepsilon'A + 4\varepsilon'^2\\
        &<(r-\varepsilon)^2+4\varepsilon'(r-\varepsilon)+4\varepsilon'^2 = (r-\varepsilon)^2+\varepsilon(r-\varepsilon)+\dfrac{\varepsilon^2}{4}\\
        &=r^2 -\varepsilon r + \dfrac{\varepsilon^2}{4}.
    \end{align*}
    Since $r>\varepsilon>\frac{\varepsilon}{4}$ then $d_{\tilde{g}_j}((x,s),(p_j,t)) < r$, that is, $(x,s) \in B_r((p_j,t))\subset N_j \times \R$. Moreover \[d_g\left( F_j(x,s),(z,s)\right) = d_g\left((f_j(x),s), (z,s)\right) = d_{S^2}(f_j(x),z) <\varepsilon' < \varepsilon.\]
\end{enumerate}
\end{proof}

\appendix

\section{Gromov-Hausdorff Convergence}\label{Appendix GH}

We briefly recall the notion of Gromov-Hausdorff convergence and some basic facts that were used throughout the paper. For a detailed exposition, we refer the reader to standard texts on Metric Geometry, for example \cite{burago2001course}.

Let $(X,d_X)$ and $(Y,d_Y)$ be compact metric spaces. The \emph{Gromov-Hausdorff distance} between $X$ and $Y$ is defined by \[d_{GH}(X,Y) \coloneqq \inf_{Z,\varphi,\psi} d_H^Z\bigl(\varphi(X),\psi(Y)\bigr)\] where the infimum is taken over all metric spaces $(Z,d_Z)$ and all isometric embeddings $\varphi\colon X\longrightarrow Z$, $\psi\colon Y\longrightarrow Z$ and $d_H^Z$ denotes the Hausdorff distance in $Z$. 

\begin{defi}
     A sequence of compact metric spaces $(X_j,d_j)$ converges in the Gromov-Hausdorff sense to a compact metric space $(X_\infty,d_\infty)$, and it is denoted by $(X_j,d_j) \xrightarrow{\ GH \ }{(X_\infty, d_\infty)}$, if $d_{GH}(X_j,X_\infty)\longrightarrow 0.$ 
\end{defi}

Equivalently, Gromov-Hausdorff convergence can be described in terms of $\varepsilon$-isometries.

\begin{defi}
    Let $X,Y$ be two metric spaces and $\varepsilon>0$. A map $f \colon X \rightarrow Y$ (possibly non continuous) is called an $\varepsilon$-isometry if its distortion \[dis(f) := \sup_{x,x' \in X} \vert d_Y(f(x),f(x')) - d_X(x,x')\vert \leq \varepsilon\] and if $f (X)$ is an $\varepsilon$-net in $Y$.
\end{defi}

\begin{prop} \label{GH isometry}
    If $d_{GH}(X,Y) < \varepsilon$, then there is a $2\varepsilon$-isometry $f \colon X \rightarrow Y$. Reciprocally if there is an $\varepsilon$-isometry $f \colon X \rightarrow Y$, then $d_{GH}(X,Y) < 2\varepsilon$.
\end{prop}

For noncompact spaces one uses the pointed version of this notion of convergence. A pointed metric space is a pair $(X,p)$ consisting of a metric space $X$ and a point $p \in X$.

\begin{defi}\label{definition pGH}
    A sequence $(X_i,p_i)_{i}$ of pointed metric spaces converges in the Gromov-Hausdorff sense to a pointed metric space $(X,p)$ if the following holds. For every $r>0$ and $\varepsilon >0$ there exists $i_0 \in \mathbb{N}$ such that for every $i>i_0$ there is a (not necessarily continuous) map $f \colon B_r(p_i) \rightarrow X$ such that:
    \begin{enumerate}
        \item[(1)] $f(p_i) = p$;
        \item[(2)] the distortion satisfies $\text{dis}\ f < \varepsilon$;
        \item[(3)] the $\varepsilon$-neighborhood of the set $f(B_r(p_i))$ contains the ball $B_{r-\varepsilon}(p)$. 
    \end{enumerate}
    In this case we write $(X_i,p_i) \xrightarrow{\ pGH \ }(X,p).$
\end{defi}

There is also the following relation between Gromov-Hausdorff and pointed Gromov-Hausdorff convergence.

\begin{prop} \label{pointed GH and GH COMPACT} For compact metric spaces the following is true:
\begin{enumerate}
    \item[(1)] $(X_n, p_n) \xrightarrow{\ pGH \ }(X,p)$ implies that $X_n \xrightarrow{\ GH \ } X$ provided there is a uniform bound on the diameter of $X_n.$
    \item[(2)] If $X_n \xrightarrow{\ GH \ } X$ and $p \in X$, then one can choose a point $p_n$ in every $X_n$ so that $(X_n,p_n)\xrightarrow{\ pGH \ }(X,p)$.
\end{enumerate}
\end{prop}

The main objects that appear as Gromov-Hausdorff limits are the so called \textit{Alexandrov spaces}. A complete length space $(X,d)$ is said to have \emph{curvature bounded below by $\kappa$ in the sense of Alexandrov} if geodesic triangles in $X$ satisfy the corresponding comparison inequalities with triangles in the simply connected $2$-dimensional space form of constant sectional curvature $\kappa$. Roughly speaking, triangles in $X$ are at least as ``fat'' as the corresponding comparison triangles in the model space. 

Every complete Riemannian manifold with bounded sectional curvature $K\geq \kappa$ is an Alexandrov space with curvature bounded below by $\kappa$. Thus the notion of Alexandrov spaces inside Gromov-Hausdorff's theory provides a natural synthetic extension of the class of Riemannian manifolds with a lower sectional curvature bound, while allowing singular spaces to arise as limits. One of the fundamental features of this class is its stability under Gromov-Hausdorff convergence.

\begin{prop}\label{limite de alexandrov} A Gromov-Hausdorff limit of Alexandrov spaces of curvature $\geq k$ is itself a space of curvature $\geq k$. The same is true for  Gromov-Hausdorff limits of pointed spaces.
\end{prop}

\begin{prop}\label{Gromov compactness}
Let $k\in \R$, $D>0$ and $n \in \mathbb{N}$. Define \[\mathfrak{M}(n,k) = \{X \colon X \text{ is an Alexandrov space of curvature} \geq k \text{ and dim}_H(X) \leq n \}\] and \[\mathfrak{M}(n,k,D) = \{X \in \mathfrak{M}(n,k) \colon \text{diam}(X) \leq D\}.\] Then the class $\mathfrak{M}(n,k,D)$ regarded with the  Gromov-Hausdorff metric is compact. In particular $\mathfrak{M}(n,k)$ is compact if $k>0$.
\end{prop}

Other properties of Alexandrov spaces are listed below, including a version of Bonnet-Myers theorem and the Bishop-Gromov theorem for these spaces under a curvature bound hypothesis.

\begin{teo}\label{Bonnet Myers Alexandrov}
    Let $X$ be an Alexandrov space of curvature $\geq k$ where $k$ is a positive number. Then $\text{diam}(X) \leq \frac{\pi}{\sqrt{k}}$.
\end{teo}

In this paper we denote $\mathcal{H}^n$ the $n$-dimensional Hausdorff measure. The Bishop-Gromov theorem in the context of Alexandrov spaces is a tool to compare the Hausdorff measure of balls in terms of their volume in spaces of constant sectional curvature.

\begin{teo}\label{Bishop Gromov alexandrov} Let $X$ be a locally compact Alexandrov space of curvature $\geq k$ and $n$ be a positive integer. Then for every $p \in X$ the ratio \[\dfrac{\mathcal H^n(B_r(p))}{V_r^k}\] is nonincreasing in $r$, where $V_r^k$ is the volume of the $n$-dimensional ball of radius $r$ in the $n$-dimensional space form $M_k^n$ of curvature $k$. In particular for $R \geq r>0$, we have: \[\dfrac{\mathcal H^n(B_R(p))}{V_R^k}\leq \dfrac{\mathcal H^n(B_r(p))}{V_r^k}.\]
\end{teo}

\begin{teo}\label{Bishop gromov equality alexandrov}
Let $n$ be a positive integer and  $X$ be an $n$-dimensional Alexandrov space of curvature $\geq k$. Then for every $p \in X$ and every $r>0$ \[\mathcal H^n(B_r(p))\leq V_r^k\] where $V_r^k$ is the volume of the $n$-dimensional ball of radius $r$ in the $n$-dimensional space form $M_k^n$ of curvature $k$. Furthermore the equality case holds if and only if $B_r(p)$ is isometric to the ball of radius $r$ in $M_k^n$. In particular, if the equality takes place for all $r>0$ then $X$ is isometric to $M_k^n$.
\end{teo}

\section{Intrinsic Flat Convergence}\label{Appendix IF}

In this section we recall the notions of intrinsic flat and volume preserving intrinsic flat convergence. This theory was introduced by Sormani and Wenger in \cite{SormaniWenger2011} in analogy to the Federer-Fleming flat distance. Ambrosio and Kirchheim in \cite{AmbrosioKirchheim2000} developed a theory of integral currents in general metric spaces, extending the classical Euclidean framework of Federer and Fleming. Rather than working with differential forms, Ambrosio and Kirchheim formulate metric currents using $(m+1)$-tuples \[\omega=(f,\pi_1,\ldots,\pi_m) \in \mathcal D^m(Z)\] where $f \colon Z\rightarrow \R$ is bounded and Lipschitz and each $\pi_i \ \colon Z\rightarrow  \R$ is a Lipschitz function defined on a complete metric space $Z$. Here $\mathcal{D}^m(Z)$, where $m$ is an integer greater or equal to one, is the set of all $(m+1)$-tuples of real valued Lipschitz functions on $Z$ whose first function is bounded. The case $m=0$ is just the set of bounded Lipschitz functions.

\begin{defi}\label{definition currents}
    An $m$-dimensional current $T$ on a complete metric space $Z$ is a real valued multilinear functional on $\mathcal D^m(Z)$ that satisfies the following conditions:
    \begin{enumerate}
        \item[i)] if there exists an index  $i \in \{1, \dots, m\}$ such that $\pi_i$ is constant on a neighborhood of $\{f \neq 0\}$, then $T(f,\pi_1,\dots,\pi_m) = 0$;
        \item[ii)]  $T$ is continuous with respect to the pointwise convergence of the functions $\pi_i$ such that its Lipschitz constant satisfies $\text{Lip}(\pi_i) \leq  1$;
        \item[iii)] there exists a finite Borel measure $\mu$ on $Z$ such that $\vert T(f,\pi_1,\dots,\pi_m)\vert \leq \prod\limits_{i=1}^m \text{Lip}(\pi_i) \int_Z \vert f\vert d\mu$ for all $(f,\pi_1, \dots, \pi_m) \in \mathcal{D}^m(Z)$.
    \end{enumerate}
   We denote by $\textbf{M}_m(Z)$ the space of $m$-dimensional currents on $Z$. The mass measure $\Vert T\Vert$ of a current $T \in \textbf{M}_m(Z)$ is the smallest Borel measure $\mu$ such that $(iii)$ holds. The mass of $T$ is then defined as:
   \[\textbf{M}(T) = \Vert T\Vert(Z) = \int_Zd\Vert T\Vert.\]
   \end{defi}

Given a Lipschitz map, it is also possible to induce a current defined by it.

\begin{defi}
    Given a Lipschitz map $\phi \colon Z \rightarrow Z'$, the pushforward of a current $T \in \textbf{M}_m(Z)$ is the $m$-current $\phi_\#T \in \textbf{M}_m(Z')$ given by:\[ \phi_\#T(f,\pi_1,\dots,\pi_m) \coloneqq T(f \circ \phi,\pi_1 \circ \phi, \dots,\pi_m \circ \phi)\] for all $(m+1)$-tuples $(f,\pi_1, \dots, \pi_m) \in \mathcal{D}^m(Z').$
\end{defi}

\begin{defi}
    The boundary of a current $T \in\textbf{M}_{m+1}(Z)$ is defined by: \[\partial T(f,\pi_1,\dots,\pi_m) \coloneqq  T(1,f,\pi_1, \dots,\pi_m)\]  Note that $\phi_\#(\partial T) = \partial(\phi_\#T)$.
\end{defi}

\begin{defi}
    Given a current $T \in \textbf{M}_{m+k}(Z)$ and a $(k+1)$-tuple $\omega = (g,\tau_1, \dots, \tau_k) \in \mathcal D^k(Z)$, the restriction $T\llcorner\omega \in \textbf{M}_{m}(Z)$ is given by \[T\llcorner\omega (f, \pi_1, \dots, \pi_m) \coloneqq T(fg, \tau_1, \dots, \tau_k, \pi_1, \dots, \pi_m).\]
\end{defi}

There are two important classes of currents: the integer rectifiable ones and the integral currents.

\begin{defi}
    A current $T \in \textbf{M}_m(Z)$ is called an integer rectifiable current, and denoted by $T \in\mathcal{I}_m(Z)$, if $T$ admits a  parametrization $(\{\phi_i\}, \{\theta_i\})_i$ such that $\{\phi_i\}$ is a countable collection of bi-Lipschitz maps $\phi_i \colon A_i \rightarrow Z$ with $A_i \subset \R^m$ a precompact Borel measurable set and pairwise disjoint images and each $\theta_i$ is a weight function $\theta_i \in L^1(A_i,\mathbb N)$ such that \[T = \sum\limits_{i=1}^\infty (\phi_i)_\#[[\theta_i]] \qquad \text{and} \qquad \textbf{M}(T) = \sum\limits_{i=1}^\infty \textbf{M}\left((\phi_i)_\#[[\theta_i]]\right)\] where $[[\theta_i]]$ is the $m$-current given by: \[[[\theta_i]](f,\pi_1, \dots, \pi_m) = \int_{A_i} \theta_i f \ d\pi_1 \wedge\cdots \wedge  d\pi_m.\] Furthermore, the mass measure of $T$ satisfies \[\Vert T\Vert = \sum\limits_{i=1}^\infty \Vert (\phi_i)_\#[[\theta_i]]\Vert.\]   
\end{defi}

\begin{defi}
An integer rectifiable current space is a triple $M=(X,d,T)$ where $(X,d)$ is a metric space and $T$ is an integer rectifiable current defined on the metric completion $\overline X$ of $X$ such that $\text{set}(T)=X$, where $\text{set}(T)$ denotes the set of points of positive lower density of the mass measure $\Vert T\Vert $. In this case it is common to denote $T = \llbracket M \rrbracket$.  
\end{defi}

\begin{defi}
    An integral current is an integer rectifiable current $T\in \mathcal I_m(Z)$ such that $\partial T$ is also a current. In this case it is denoted $T \in \textbf{I}_m(Z)$. 
\end{defi}

\begin{rmk}
     Given an integral current $T \in \textbf{I}(Z)$, by Theorem 8.6 in \cite{AmbrosioKirchheim2000}, $\partial T$ is also an integral current.
\end{rmk}

\begin{defi}
    An $m$-dimensional integral current space is an integer rectifiable current space $(X,d,T)$ whose current structure, $T$, is an integral current, that is, $T \in \textbf{I}_m(\bar{X})$. If $\partial T =0$ we say that $(X,d,T)$ is an integral current without boundary. We denote $\mathcal M^m$ the space of $m$-dimensional integral current spaces.
\end{defi}

Every oriented compact Riemannian manifold $(M^n,g)$ with boundary determines an integral current space $(M,d_g,T_M)$ where $d_g$ is the intrinsic Riemannian distance and $T_M$ is the current given by integration over $M$ with respect to the chosen orientation. Its mass agrees with the Riemannian volume $\mathbf M(T_M)=\text{Vol}_g(M)$ and the boundary $\partial M$ defined as a current space agrees with the definition of $\partial M$ in Riemannian geometry.

Once all these notions have been established, it is possible to defined, given two integral current spaces, the intrinsic flat distance between them. More precisely, let $M_1, M_2 \in \mathcal{M}^n$, that is, $M_1 = (X_1,d_1,T_1)$ and $M_2=(X_2,d_2,T_2)$ are $n$-dimensional integral current spaces. Their \emph{intrinsic flat distance} is defined by \[d_{\mathcal F}(M_1,M_2) = \inf_{Z,\varphi_1,\varphi_2} d_F^Z\left(\varphi_{1\#}T_1,\varphi_{2\#}T_2\right) \] where the infimum is taken over all complete metric spaces $(Z,d_Z)$ and all isometric embeddings $\varphi_i \ \colon \ (\bar{X}_i,d_i)\longrightarrow (Z, d_Z)$. Here $d_F^Z$ denotes the flat distance between currents in the common ambient space $Z$. Explicitly, \[d_F^Z(S_1,S_2) = \inf \{\textbf{M}(U) + \textbf{M}(V)\  \colon \ U \in \textbf{I}_n(Z),\ V \in \textbf{I}_{n+1}(Z), \ S_1-S_2 = U+\partial V  \}.\]

\begin{defi}
     A sequence of integral current spaces $M_j$ converges to $M_\infty$ in the \emph{intrinsic flat sense}, and in this case we write $M_j \xrightarrow{\mathcal F}M_\infty$, if $d_{\mathcal F}(M_j,M_\infty)\xrightarrow{\ j \rightarrow +\infty \ } 0.$
\end{defi}

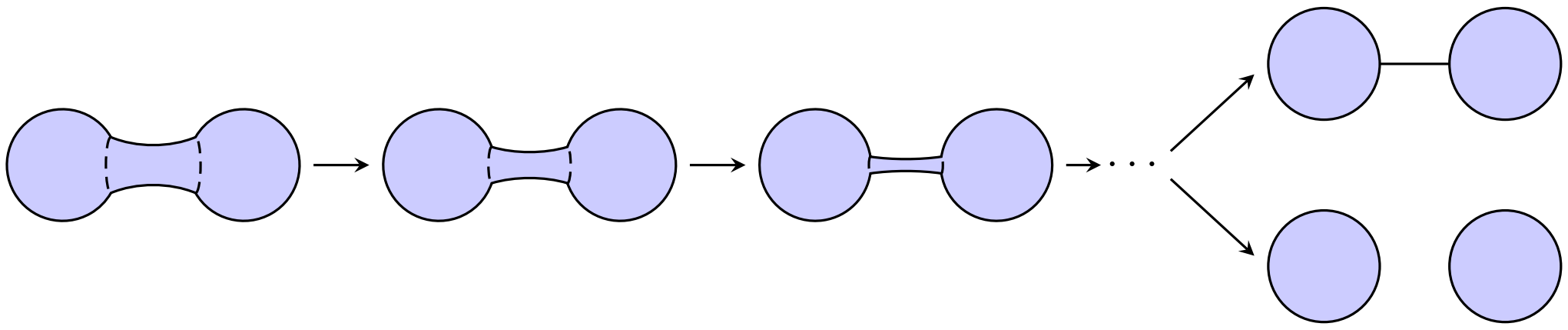
\begin{figure}[ht]
  \centering
\begin{tikzpicture}[scale=0.6, line width=0.6pt, >=stealth]
  \SphereNeckSphere{0}{0.40}
  \draw[->] (3.6,0) -- (4.4,0);
  \SphereNeckSphere{5.4}{0.26}
  \draw[->] (9.0,0) -- (9.8,0);
  \SphereNeckSphere{10.8}{0.12}
  \draw[->] (14.4,0) -- (14.9,0);
  \node at (15.4,0) {$\cdots$};

  \draw[->] (15.9,0.2)  -- (17.1,1.3);
  \draw[->] (15.9,-0.2) -- (17.1,-1.3);

  \SpheresSegment{18.1}{1.45}
  \SpheresApart{18.1}{-1.45}

\end{tikzpicture}
  \caption{A sequence of manifolds with distinct Gromov-Hausdorff (upper right) and intrinsic flat limits (lower right).}
    \label{GH vs IF}
\end{figure}

The intrinsic flat topology is well adapted to situations in which thin regions have small volume. In contrast with Gromov-Hausdorff convergence, portions of the space whose mass tends to zero may disappear completely in the intrinsic flat limit. Moreover, because the convergence takes into account the current structure and orientation, cancellation phenomena may occur. Thus the intrinsic flat limit need not coincide with the Gromov-Hausdorff limit, even when both limits exist (see Figure \ref{GH vs IF} for an illustrative example). 

An important property of integral currents is the semicontinuity: if $T_j$ converges weakly to $T_{\infty}$, then: \[\mathbf M(T_\infty) \leq \liminf_{j\to\infty}\mathbf M(T_j).\] For Riemannian manifolds, this property becomes \[\text{Vol}(M_\infty)\leq\liminf_{j\to\infty}\text{Vol}(M_j).\] In particular, intrinsic flat convergence alone allows loss of volume in the limit. To rule out such loss, one considers the stronger notion of \emph{volume preserving intrinsic flat convergence}.

\begin{defi}
  We say that a sequence of integral current spaces $M_j=(X_j,d_j,T_j)$ converges to $M_\infty\coloneqq(X_\infty,d_\infty,T_\infty)$ in the volume preserving intrinsic flat sense if \[M_j\xrightarrow{\mathcal F}M_\infty\qquad \text{and} \qquad \mathbf M(T_j)\longrightarrow \mathbf M(T_\infty).\] In this case we write $M_j\xrightarrow{\mathcal{VF}}M_\infty.$
\end{defi}

For oriented Riemannian manifolds, the volume preserving intrinsic flat convergence is equivalent to requiring \[d_{\mathcal F}(M_j,M_\infty) \longrightarrow 0 \qquad \text{together with} \qquad \text{Vol}_{g_j}(M_j)\longrightarrow\text{Vol}(M_\infty).\]

Integral currents also have a compactness theorem due to Wenger \cite{Wenger2011}. Roughly speaking, if a sequence of $n$-dimensional integral current spaces satisfies uniform bounds on their masses and also for their boundaries together with a uniform diameter bound, then there exists a subsequence converging in the intrinsic flat sense to an integral current space, possibly the zero space. In the case of closed oriented Riemannian manifolds, the boundary term vanishes so uniform volume and diameter bounds yield intrinsic flat subsequential compactness.
\begin{teo}[Theorem 1.2 in \cite{Wenger2011}]\label{Wenger compactness theorem} Let $m,C,D>0$ be positive constants and $M_j = (X_j, d_j, T_j)$ a sequence of integral current spaces such that each $T_j \in \textbf{I}_m(\overline{X}_j)$ satisfies \[\textbf{M}(T_j) + \textbf{M}(\partial T_j) \leq C \qquad \text{and} \qquad \text{diam}(\text{spt}(T_j)) \leq D,\] 
where $\text{spt}(T_j)$ is the support of the current $T_j$ which is defined as
\[\text{spt}(T_j) \coloneqq \text{spt}\Vert T_j\Vert = \{p \in \bar{X}_j \colon \Vert T_j \Vert (B_r(p)) > 0 \text{ for all } r > 0\}.\]
Then there exists a subsequence which converges in the intrinsic flat sense to an integral current space of the same dimension. 
\end{teo}

\vspace{0.25cm}

\bibliographystyle{plain}
\bibliography{bibliography.bib}

\end{document}